\documentclass[12pt,twoside,leqno]{amsart}
\usepackage{amsmath,amscd,amssymb,amsfonts}

\newcommand{\N}{\mathbb{N}}
\newcommand{\Z}{\mathbb{Z}}
\newcommand{\Q}{\mathbb{Q}}

\newcommand{\C}{\mathbb{C}}
\renewcommand{\H}{\mathbb{H}}
\renewcommand{\P}{\mathbb{P}}
\renewcommand{\L}{\mathbb{L}}

\newcommand{\mE}{\mathcal{E}}
\newcommand{\mF}{\mathcal{F}}

\newcommand{\mO}{\mathcal{O}}
\newcommand{\mP}{\mathcal{P}}
\newcommand{\mY}{\mathcal{Y}}
\newcommand{\mZ}{\mathcal{Z}}

\newcommand{\lra}{\longrightarrow}

\theoremstyle{plain}

\numberwithin{thm}{section}
\numberwithin{equation}{section}

\begin{document}
\title{Calabi--Yau Operators\\\hfill\\ }
\author{Duco van Straten}
\date{}
\thanks{This work is entirely supported by DFG Sonderforschungsbereich/Transregio 45}
                         %= Fussnotentext.
\begin{abstract}
Motivated by mirror symmetry of one-parameter models, an interesting
class of Fuchsian differential operators can be singled out, the so-called
{\em Calabi--Yau operators}, introduced by {\sc Almkvist} and {\sc Zudilin} in \cite{AZ}. They conjecturally determine $Sp(4)$-local systems that underly a
$\Q$-VHS with Hodge numbers
\[h^{3 0}=h^{2 1}=h^{1 2}=h^{0 3}=1\]
and in the best cases they make their appearance as Picard--Fuchs
operators of families of Calabi--Yau threefolds with $h^{12}=1$ and encode the numbers of rational curves on a mirror manifold with $h^{11}=1$.
We review some of the striking properties of this rich class of operators.

\end{abstract}

\maketitle
\section{Calabi--Yau operators}

\centerline{\bf \em The story of the quintic}
\vskip 10pt
The story of  Calabi--Yau operators is connected to the beginnings of
mirror symmetry, in particular with the classical paper by {\sc Candelas},
{\sc de la Ossa}, {\sc Green} and {\sc Parkes} \cite{COGP}, which is still
an excellent introduction to the subject.
The larger story how mirror symmetry  entered the mathematical community and has shaped a good part of present day mathematics has been told in more detail at other places, and we refer to \cite{Mor2,Tho, Yau2} for nice surveys and \cite{Katz2, Voi, CK, Horics} for a more comprehensive accounts of this ever growing subject. In this paper we can only give the barest outline as far as relevant for our purpose.\\

Let us start with recalling the mysterious calculation with
the power series
\[ y_0(t)=\sum_{n=0}^{\infty} \frac{(5n)!}{(n!)^5} t^n =1+120 t+113400 t^2+\ldots \in \Z[[t]] \]
that appeared in \cite{COGP}.
It represents the unique (normalised) holomorphic solution to
the hypergeometric differential operator
\[ \mathcal{P}:=\Theta^4-5^5 t (\Theta+\frac{1}{5})(\Theta+\frac{2}{5})(\Theta+\frac{3}{5})(\Theta+\frac{4}{5}) \in \Q \left[t,\frac{d}{dt}\right], \]
where
\[ \Theta:=t\frac{d}{dt}\]
denotes the logarithmic derivative with respect to the parameter $t$.
By expressing the operator in a new coordinate $q$, we can bring $\mathcal{P}$ to a normal form
\[ \mathcal{P}=\theta^2 \frac{5}{K(q)} \theta^2,\;\;\; \theta:=q \frac{d}{dq},\]
where $K(q)$ is a power series in $q$.  In fact, it is easy to see that this $q$-coordinate is given by
\[ q = e^{y_1(t)/y_0(t)}=t+770t^2+\ldots,\]
where
\[y_1(t):=\log(t) y_0(t)+f_1(t),\;\;f_1(t) \in t\hspace{1pt}\Q[[t]] \]
is the normalised solution of $\mathcal{P}$ that contains a single logarithm.
It is easy to compute the beginning of the power series expansion of $K(q)$ and write it in the form
of a Lambert series
\[K(q)=5+\sum_{d=1}^{\infty} n_d \frac{d^3q^d}{1-q^d}\]
from which one can read off numbers $n_d$, which a priori are in $\Q$.
One finds
\[ n_1=2875,\;\;\;n_2=609250,\;\;\;n_3=317206375,\;\;\ldots\]
and computing the $n_d$'s further, it appears that these numbers are all integers.\\

What makes the above calculation intriguing is the fact that it is related to the properties of two very different {\em Calabi--Yau threefolds}, but which the physics of strings suggests to be closely related.\\

{\bf A-incarnation}: The first manifold is the general quintic hypersurface
$X \subset\P^5$, which is a Calabi--Yau threefold with Hodge numbers $h^{11}=1$ and $h^{12}=101$.
The numbers $n_d$ are called {\em instanton numbers} and were argued in \cite{COGP} to be equal to the {\em number of rational
degree $d$ curves on $X$} counted in an appropriate way. This was a big claim, as only the first two
numbers were known at the time: the number $2875$ of lines on a general quintic
was determined by the founding father of enumerative geometry {\sc H. Schubert} in 1886 \cite{Schu}, the number of $609250$
of conics was determined by {\sc S. Katz} \cite{Katz1} hundred years later. In a heroic {\em tour the force}, the number of twisted cubics on the quintic was determined by {\sc S. Str{\o}mme} and {\sc G. Ellingsrud} \cite{ES} and in fact served as a crucial cross-check for the above calculation and resulted in the famous message: {\em Physics wins}! For the details of that story we refer to \cite{Xam}.\\

{\bf B-incarnation}: The second manifold is the {\em quintic mirror $Y$}
with 'flipped' Hodge numbers $h^{12}=1, h^{11}=101$. It was constructed via
an orbifold construction that was proposed earlier by {\sc B. Greene} and {\sc R. Plesser} \cite{GP}. As the Hodge number $h^{12}(Y)$ is equal to the dimension of the local moduli space of $Y$, we have in fact a $1$-parameter family of
manifolds $Y_t$, parametrised by $t$. It can be obtained from the quintics of the so-called {\em Dwork pencil}
\[ \sum_{i=1}^5 x_i^5-5 \psi \prod_{i=1}^5 x_i=0,\;\;\;t=\frac{1}{(5 \psi)^5} \]
by dividing out the abelian group of order $125$ generated by $x_i \mapsto \zeta_i x_i, \zeta_i^5=1, \prod_{i=1}^5 \zeta_i =1$ and resolving the resulting
singularities. The solution to the differential equation $y_0(t^5)$ is in fact a (normalised) {\em period integral} of $Y_t$ and $\mathcal{P}$ is the associated {\em Picard--Fuchs equation}.
The most salient property of the differential operator is the fact that it has a so-called MUM (=maximal unipotent monodromy) point at $0$, where the
variety $Y_t$ degenerates to a union of divisors forming a combinatorial sphere.\\

The computation of {\sc Candelas} and coauthors was immediately extended to other Calabi--Yau threefolds,\cite{Mor1, Fon, KT1, KT2, LT}. The extension to smooth complete intersections in weighted projective spaces yield $13$ cases for which the associated differential equation is of hypergeometric type, \cite{BvS}. In fact, there is a $14$th case \cite{Alm1, DoMo} 
that was not considered at that time, as it corresponds to a Calabi--Yau variety with a singular point. \\

In $1993$ {\sc Victor Batyrev} came up with a general interpretation of mirror symmetry in terms of dual reflexive polytopes \cite{Bat1}, which lead to a plethora of examples of mirror pairs of Calabi--Yau manifolds, but usually with rather large Hodge
numbers. As the mirror manifold $Y$ of $X$ varies in a $h^{11}(X)=h^{12}(Y)$-dimensional
family, the period integrals are solutions to differential systems with this
many variables. In general only in the cases that $h^{11}(X)=h^{12}(Y)=1$, the so-called
{\em one-parameter models}, one obtains a single ordinary differential
equation of order four annihilating the period integrals. However, in \cite{BvS} first examples of Calabi--Yau threefolds with $h^{11}(X) >1$ were considered, which, due to a symmetry, still lead to a fourth order operator.  A nice example is the case of the degree $(3,3)$ hypersurface $X$ in $\P^2 \times \P^2$, which leads to the operator number $15$ in the AESZ-list \cite{AESZ}
\[\Theta^4  - 3 t (3\Theta + 1)(3\Theta + 2)(7\Theta^2+ 7\Theta + 2)
  -  72 t^2(3\Theta + 5)(3\Theta + 4)(3\Theta + 2)(3\Theta + 1) ,\]
which no longer is of hypergeometric type.  A curve on $X$ has two degrees, coming from the two $\P^2$-factors. The corresponding instanton numbers of the above operator count the rational curves with total degree equal to $d$.
The mirror manifold of $X$ has $h^{12}=2$, but over a line in the two-dimensional deformation space the cohomology splits off a sub Hodge structure with $h^{03}=h^{12}=h^{21}=h^{30}=1$.\\

The discovery that the computation of periods of one manifold provides
enumerative information about another manifold was totally unexpected and
left people wonder about the geometrical relation between $X$ and $Y$.
It was a key motivation for the development of mathematical understanding of
mirror symmetry and led to several important insights. First, {\em Gromov--Witten theory} was developed to provide a rigorous basis for counting curves on general manifolds, \cite{KM}. This enabled {\sc Givental} \cite{Giv1} and {\sc Lian, Liu} and {\sc Yau} \cite{LLY} to prove the {\em mirror theorem} that vindicated the above computational scheme, but left out the question of the geometrical relation between the spaces $X$ and $Y$. Second, the idea that in mirror symmetry the {\em symplectic geometry} of $X$ gets identified with the
{\em holomorphic geometry} of $Y$ and vice versa got a precise expression in terms of {\sc Kontsevich}'s notion of {\em homological mirror symmetry}, \cite{Kon}. The insight that this in turn leads to a description of $X$ and $Y$ as two dual
torus fibrations by {\sc Strominger--Yau--Zaslow} \cite{SYZ} took some of the mystery of the  mirror symmetry phenomenon, but left the mathematical community with very difficult problems to solve. The approach of {\sc M. Gross} and {\sc B. Siebert} seeks to develop this picture of mirror duality in the framework of algebraic geometry out of dual logarithmic degeneration data and the resulting affine manifolds with singularities, \cite{GS1, GS2}, which can be seen as a grand generalisation of {\sc Batyrev}'s notion of dual reflexive polytopes.\\

\centerline{\bf \em Quantum Cohomology at the Mittag-Leffler Institute 1996--1997}
\vskip 10pt
In the year 1996/97 a special year on {\em Enumerative Geometry and its Interaction with Theoretical Physics} was organised by Geir Ellingsrud, Dan Laskov, Anders Thorup and Stein Arild Str{\o}mme at the Mittag-Leffler Institute.
The text  \cite{Alu} collects write-ups of the talks that were given during the first semester and capture very well the exciting atmosphere aroused by the new techniques of {\em Gromov--Witten theory}, {\em Frobenius manifolds}, {\em quantum cohomology}, {\em quantum $D$-modules}, all aimed at understanding the mirror theorem.\\

During my stay at Mittag-Leffler, I intended to steer away from toric mirror
symmetry and tried to obtain further examples of one-parameter models by
looking at complete intersections in other simple spaces, like homogeneous
spaces.  For example, complete intersections in Grassmanianns lead to varieties with $h^{11}=1$, but as these were not of toric type, it was not so clear how to
obtain a mirror dual, let alone its Picard--Fuchs equation.
This had changed after {\sc Dubrovin} \cite{Dub} and {\sc Givental} \cite{Giv2} showed that it is
possible to find the Picard--Fuchs equation directly from $X$ in terms of Gromov--Witten invariants of an ambient manifold of $X$. More precisely, if $X$ is a complete
intersection in a manifold $Z$ that is simple enough to allow for an explicit
description of its quantum cohomology ring, one can use a {\em Laplace transform} to obtain the differential equation for $Y$: the {\em quantum Lefschetz principle} in the formulation of \cite{Giv2, Giv3}.\\

To explain these important ideas, consider for simplicity a smooth projective variety $Z$ with $h^2(Z)=1$,  without odd cohomology  and let $H \in H^2(Z)$ be its ample generator.
The homology class of a curve  $\Sigma$ in $Z$ is determined by its
{\em degree}, defined as the intersection number $\Sigma\cdot H$.
For $A,B,C \in H^*(Z)$ the {\em Gromov--Witten three-point function} is
the series
\[ \langle A,B,C \rangle:=\sum_{d=0}^{\infty} \langle A,B,C \rangle_dt^d ,\]
where
\[ \langle A,B,C \rangle_d\]
is the Gromov--Witten count of rational degree $d$ curves that meet the
cycles (Poincar\'e dual to) $A,B$ and $C$. The quantum product $\star$ is
the $t$-dependent product determined by the equation
\[ (A\star B,C) =\langle A,B,C \rangle ,\]
where on the left we use the non-degenerate Poincar\'e pairing on $H^*(Z)$.
The product $\star$ is associative and commutative and in the Fano case
can be used to  define a new ring structure on $QH^*(Z):=H^*(Z)[t]$:
the quantum cohomology ring of $Z$.
The {\em Dubrovin--Givental} connection is the connection
\[ \nabla=td-H\star \]
on the trivial bundle over the $t$-line and with the vector space $H^*(Z)$
as fibre, whose horizontal sections $S(t)$ are solutions to the differential
system
\[ \theta S=H \star S,\;\;\theta:=t \frac{d}{dt} .\]

For $Z=\P^4$ the only non-trivial three-point invariant is
\[ \langle H^4, H^4, H \rangle_1=1 ,\]
which expresses the obvious enumerative fact that there is a single line
through two points and this line intersects a given hyperplane in a single
point. The quantum cohomology ring is identified with $\C[H,t]/(H^5-t)$, i.e. 
one has
\[ H \star H=H^2,\;\; H^2 \star H=H^3,\;\; H^3 \star H=H^4,\;\; H^4\star H=t\]
Using the basis $1, H, H^2, H^3,H^4$ for $H^*(\P^4)$, the quantum differential system for $S=\sum_i S_i H^i$ can be written as
\[ t \frac{d}{d t} \left( \begin{array}{c}
S_0\\S_1\\S_2\\S_3\\S_4 \end{array} \right)=
\left( \begin{array}{ccccc}
0&0&0&t\\
1&0&0&0\\
0&1&0&0\\
0&0&1&0\\
0&0&0&1
\end{array} \right) \cdot
\left( \begin{array}{c}
S_0\\S_1\\S_2\\S_3\\S_4 \end{array}\right),
\]
which leads to the differential equation of order $5$ for the lowest component $S_4(t)$ of $S(t)$:
\[ (\theta^5-t)S_4(t)=0 .\]
This equation is easily seen to have
\[ \psi(t):=\sum \frac{1}{(n!)^5}t^n\]
as unique holomorphic solution.
The {Laplace transform} of this function
\[\frac{1}{t}\int_0^{\infty} \psi(s^5) e^{-s/t}ds\]
is the function
\[ \phi(t^5) =\sum_{n=0}^{\infty} \frac{(5n)!}{(n!)^5}t^{5n}=y_0(t^5) .\]
Note that the $5$'s are dictated by the fact that the canonical class of
$\P^4$ is $K=-5H$ and it transforms the irregular quantum differential
system into one with regular singularities!\\

A closely related aspect was the idea, already present in \cite{CV}, that the notion of mirror symmetry should be extended from Calabi--Yau spaces to the ambient Fano manifolds and that these mirrors were described by a {\em Landau--Ginzburg potential}, for example in case of $\P^4$ by the Laurent polynomial
\[ W=X_1+X_2+X_3+X_4 + \frac{1}{X_1X_2X_3X_4} . \]
The solutions to the quantum differential system of $\P^4$ have a
representation as oscillatory integrals attached to $W$, which by
Laplace transformation become period integrals of the manifold
\[  \{ 1-t W=0 \}  \subset (\C^*)^4\]
that completes to a Calabi--Yau space by compactification in the
toric manifold defined by the Newton polytope of $W$.
The upshot is the following: by expanding into a geometric
series, the normalised period
\[\frac{1}{(2\pi i)^4}\oint \frac{1}{1-tW}\frac{}{}\frac{dX_1}{X_1}\frac{dX_2}{X_2}\frac{dX_3}{X_3}\frac{dX_4}{X_4}\]
expands as
\[ \sum_{n=0}^{\infty} [W^n]_{0} t^n, \]
where $[-]_0$ takes the constant term of a Laurent polynomial.
Indeed, for the above Laurent polynomial for the mirror of $\P^4$ one
obtains
\[ [W^{5n}]_{0} = \frac{(5n)!}{(n!)^5} ,\]
which leads back to our function $y_0(t^5)$. For a much more detailed analysis
of this example see \cite{Giv4, Gue1, Gue2}.\\

Now for the Grassmannian $Z=G(n,k)$ the quantum cohomology was determined
by {\sc Siebert--Tian}, \cite{ST} and during the special year at Mittag-Leffler the idea arose to use this to calculate the Laplace transform of
the quantum differential  system of the Grassmannian and try to come up with 
predictions for the number of rational curves on the Grassmannian Calabi--Yaus. This exciting collaboration with {\sc Batyrev}, {\sc Kim} and {\sc Ciocan-Fontanine} led to the papers \cite{BCKvS1} and \cite{BCKvS2}.\\

As an example, consider the Grassmannian $Z=G(2,5)$; it is a six-dimensional
Fano variety with $K=-5H$. From the quantum cohomology differential system
we obtained in \cite{BCKvS1} the function
\[\psi(t):=\sum_{n=0}^{\infty} \frac{1}{n!^5} A_n t^n ,\]
where
\[ A_n:=\sum_{k=0}^n {n \choose k}^2 {n+k \choose k}\]
are the famous Ap\'ery numbers related to $\zeta(2)$, \cite{Ape}.
The complete intersection $X:=X(1,2,2) \subset Z$ of
$Z$ with three general hypersurfaces of degree $1, 2, 2$ in the Pl\"ucker
embedding of the Grassmannian is a Calabi--Yau threefold with $h^{11}=1$.
Using the quantum Lefschetz principle/Laplace transform we were lead to
the function
\[\phi(t)=\sum_{n=0}^{\infty} \frac{n! (2n)! (2n)!}{(n!)^5} A_n t^n=\sum_{n=0}^{\infty} {2n \choose n}^2 A_n t^n\]
that should be a normalised period of a mirror $Y$ to $X(1,2,2)$.
The function $\phi$ is the holomorphic solution to the operator
\[ \mathcal{P}=\Theta^4-4t(2\Theta+1)^2(11\Theta^2+11\Theta+3)-16t^2(2\Theta+1)^2(2\Theta+3)^2\]
with Riemann symbol
\[
\left\{
\begin{array}{cccc}
0&\alpha&\beta&\infty\\
\hline
0&0&0&1/2\\
0&1&1&1/2\\
0&1&1&3/2\\
0&2&2&3/2
\end{array}
\right\} ,
\]
which is number 25 in the AESZ-list \cite{AESZ}.
The instanton numbers of $X(1,2,2)$ were then found to be
\[n_1=400, n_2=5540, n_3=164400,\;\;n_4=7059880,\ldots \]

Of particular interest is the case of the Calabi--Yau section
$X:=X(1,1,1,1,1,1,1) \subset G(2,7)$ obtained by taking $7$ generic
linear sections of the Grassmannian.
We found that in this case the function $\psi$ is given as
\[\psi(t)=\sum_{n=0}^{\infty} A_n t^n,\]
where
\[ A_n=\sum_{k,l} {n \choose k}^2 {n \choose l}^2{k+l \choose n}{2n-k \choose n} . \]

At the same time at Mittag-Leffler, {\sc Einar R{\o}dland} determined the
mirror for the generic Pfaffian Calabi--Yau in $ X' \subset \P^7$ via the
orbifold method \cite{Roe}. When he showed me the operator he had obtained,
we were both electrified: it was {\em identical} with the above 
$G(2,7)$-operator that I had obtained a week before.
But the instanton numbers for $X$ and $X'$ had to be different! The mystery
was partly resolved after realising that the operator had
{\em two points of maximal unipotent monodromy}.
At the origin we get the instanton numbers for for the Grassmannian Calabi--Yau $X$,
\[n_1 = 196,\;\; n_2 = 1225,\;\; n_3 = 12740,\;\; n_4 = 198058,\;\; n_5 = 3716944,\;\;\ldots\]
and for the point at infinity the instanton numbers for the Pfaffian Calabi--Yau
$X'$:\\
\[n_1 = 588,\;\; n_2 = 12103,\;\; n_3 = 583884,\;\; n_4 =41359136 ,\;\; n_5 =360939409,\;\;\ldots\]
At the time we were left to wonder about the geometrical relation between $X \subset G(2,7)$ and the
Pfaffian $X' \subset \P^7$. They are not birational, but it was shown later 
in \cite{Bor} that the derived categories of $X$ and $X'$ are equivalent, as 
predicted by homological mirror symmetry.\\

Although these examples are not toric, it turns out that mirror
symmetry for these examples still can be linked up with {\sc Batyrev}'s
theory of dual reflexive polyhedra: the Grassmanianns can be
{\em degenerated} to a toric variety with singularities in codimension 3,
\cite{BCKvS1}.
This leads to a Laurent polynomial description for the Grassmannian that was
found before by {\sc Eguchi}, {\sc Hori} and {\sc Xiong}, \cite{EHX}. For
a beautiful recent approach to the mirror symmetry of the Grassmannian, its
relation to the Langlands dual group and the cluster structure, 
see \cite{MaRi} and \cite{RiWi}.\\

It was suggested at the time that one could try to invert the degeneration
construction and start with special singular toric manifolds and smooth these
to obtain further examples, see \cite{Bat2} and \cite{AS}. It has been verified that all Fano varities of dimension $2$ and $3$ admit such toric degenerations.\\

\centerline{\bf \em Calabi--Yau operators}
\vskip 5pt
In $2003$ I received a letter from {\sc Gert Almkvist} in which he asked if I knew more
operators {\em like the one for the quintic}. Apart from the cases coming from \cite{BvS} and the Grassmannian cases from \cite{BCKvS1}, I knew a few more coming from the construction in \cite{AS}, but soon ran out of further examples.
Then, by insightful playing with various sums of binomial coefficients, {\sc Almkvist} discovered
many further examples. In the paper \cite{AZ} of {\sc Almkvist} and {\sc Zudilin} the notion of
{\em Calabi--Yau operator} was formulated, which is more or less characterised by the condition that the calculation of \cite{COGP} works. The operators were collected in a list \cite{AESZ}.
For a slightly more systematic and updated list, see \cite{AvS} and the online
database \cite{CYDB}. \\

{\bf Preliminary definition:}
An irreducible fourth order differential operator $\mathcal{P}\in \C\left[t,\frac{d}{dt}\right]$ is called a {\em Calabi--Yau operator} if
it satisfies the following conditions:
\begin{itemize}
\item it is of Fuchsian type.
\item it is self-dual.
\item it has $0$ as MUM-point.
\item it satisfies integrality conditions:
\begin{itemize}
\item the holomorphic solution $y_0(x) \in \Z[[x]]$.
\item the $q$-coordinate $q(x) \in \Z[[x]]$.
\item the instanton numbers $n_d \in \Z$.
\end{itemize}
\end{itemize}

In fact, it is more natural to allow mild denominators and look for $N$-integral solution, $q$-coordinate and instanton numbers. Also, one could replace $\Z$ by rings of integers in a number field. There is a natural notion of Calabi--Yau operator of arbitrary order that we will not spell out here. For a more thorough discussion we refer to the thesis of {\sc Bogner} \cite{Bog} and \cite{Bog2}.
Operators of order two tend to come from families of  elliptic curves, those of order three are obtained, by a famous theorem of {\sc Fano} \cite{Fan}, from those of order two by taking the second symmetric power and appear as Picard--Fuchs operators for  families of K3-surfaces with Picard number equal to $19$. So operators of order two and three belong to more classical realms of algebraic geometry and modular forms. The case of fourth order operators seems to be the first
that leads us into completely unknown territory.

It follows from the self-duality condition of $\mathcal{P}$ that
there is a unique formal coordinate transformation $x \mapsto q=x+\ldots$ called
the {\em mirror map} that brings the operator in the form
\[ \theta^2 \frac{1}{K(q)}\theta^2,\]
where $K(q)=1+\ldots$ is a power series, called the normalised Yukawa coupling of $\mathcal{P}$.
This power series is an invariant of the operator, unique up to a scaling in $q$.
The (normalised) instanton numbers
\[ n_1,\;\;n_2,\;\;n_3,\ldots \in \Q\]
of the operator $\mathcal{P}$ are defined by writing the Yukawa coupling in the form
\[K(q)=1+\sum_{d=1}^{\infty} n_d d^3 \frac{q^d}{1-q^d} .\]

A general construction to obtain a self-dual fourth order operator is  by taking the symmetric cube of a second order operator. But these are very special as for these operators all instanton numbers $n_d$
vanish for $d>1$ and are thus counted as {\em trivial}. (Although they play a role in mirror symmetry for abelian varieties, see \cite{BvS} for an example.)
The conditions are not all independent of each other; for example the integrality properties already {\em imply} the Fuchsian nature of the differential equation, \cite{And}.\\

There is a couple of obvious questions one can ask:\\

{\bf \em Question 1:} How to construct examples of Calabi--Yau operators?\\

This question is partly but eloquently answered in the paper 
{\em The art of finding Calabi--Yau differential equations} 
by {\sc G. Almkvist}, see \cite{Alm2}.\\

We summarise here the basics. It is rather easy to fulfill the first three conditions, but
the integrality conditions are much harder to satisfy.\\

{\em A-incarnations}: As explained above, by looking at the quantum cohomology
of a Fano manifold $Z$ one obtains a quantum differential system. If $Z$ contains a Calabi--Yau threefold (and the Picard number is one), one  obtains a
Calabi--Yau operator by Laplace transformation. As we do not know the complete
classification of Fano manifolds in dimension $\ge 4$ one can not go far
beyond a class of obvious examples, obtained from homogeneous bundles over
homogeneous spaces. In fact, in \cite{EvS} it was suggested that by reverse
engineering one could make predictions about the existence of manifolds with
given characteristic numbers from the monodromy of the operator alone.\\

{\em B-incarnations}: Here one is in much better shape. Take any interesting
looking family of Calabi--Yau threefolds with $h^{12}=1$ and compute the
Picard--Fuchs operator. The chances are good that it will have a MUM-point
somewhere. Below we describe two classes of examples we have been
looking at recently. This approach is far from exhausted and there are
many more constructions one can try. The recent algorithm of {\sc Lairez} \cite{Lai} for finding Picard--Fuchs operators is most useful here.\\

{\em Hadamard products}: If $f(t)=\sum_n a_n t^n$ and $g(t)=\sum b_n t^n$ are
power series, the series
\[ f\star g (t):=\sum a_n b_n t^n\]
obtained by taking the coefficentwise product is called the {\em Hadamard product}. A classical theorem of {\sc Hadamard} states that if $f$ and $g$ satisfy a Fuchsian differential equation, then so does $f \star g$. In this way quite a few Calabi--Yau operators were found. On the level of local systems, this comes down to taking the (multiplicative) {\em convolution} of the corresponding local systems. By the
work of {\sc Katz}, {\sc Dettweiler}, {\sc Reiter}, {\sc Sabbah} \cite{NKatz, DR, DS} the monodromy and Hodge numbers of such convolutions are
under explicit control. This answers completely a question posed at the end of
\cite{EvS}.\\

{\em Binomial Sums:} In many examples the coefficients $a_n$ of the
holomorphic solution are given as special binomial sums. {\sc Almkvist} is the
uncontested champion in guessing binomial sums that give Calabi--Yau operators.
Using {\tt Zeilberger} in  {\sc Maple} allows one to find the recursion and
hence the Picard--Fuchs operator effectively.\\

{\em Computer search:} One can start with a parametric differential
equation and make a computer search for those which give rise
to cases with integral solution, mirror map and instanton numbers.
This was done in \cite{AESZ} for operators of degree two. Going to
higher degree might be possible, but is hampered by the fact that
the number of free parameters becomes too big to handle by brute force.\\

{\em Pullback from fifth order:} A characteristic property of Calabi--Yau operators is the vanishing of a certain quantity $Q$ (see section 2.4), which causes the second order Wronskians
to satify a fifth order operator with MUM at the origin. The fourth and fifth order operator determine each other; on the level of Lie algebras this is the exceptional isomorphism  $sp(4) \approx so(5)$. One starts from fifth order operators and finds by ``pullback'' the corresponding fourth order operator. A couple of Calabi--Yau operators were found this way, but there seem to be very few simple fifth order operators that can be used.\\

The relation between fourth and fifth order operators was used with
great success by {\sc Almkvist} and {\sc Guillera} \cite{AG, ABG} to find
Ramanujan type formulas for $\frac{1}{\pi^2}$, which is a formula of type
\[
\sum_{n=0}^{\infty }A_{n} (a+bn+cn^{2}) z_{0}^{n}=\frac{1}{\pi ^{2} .}
\]
The first were found by {\sc Guillera} and five of them were proved by using the {\sc Wilf--Zeilberger} machinery. It was later realized that $A_{n}$ was the coefficient of hypergeometric fifth order Calabi--Yau equation. Later, using the properties of the fourth degree pullback, {\sc Almkvist} and {\sc Guillera} found several more formulas also for non-hypergeometric equations. A striking example is the formula%
\[
\sum_{n=0}^{\infty }\frac{(6n)!}{n!^{6}}(36+504n+2128n^{2})\frac{1}{%
1000000^{n}}=\frac{375}{\pi ^{2}}
\]%
which in principle can be used to compute an arbitrary decimal of $\dfrac{1}{\pi ^{2}}$
without computing the earlier ones.\\

{\em Laurent series}: This is a  special case of a B-incarnation.
From a Laurent polynomial $W$ we can compute the constant term series
\[ \sum_{n=0}^{\infty} [W^n]_{0}t^n \]
and from it one can in turn find the Picard--Fuchs operator that annihilates it.
(In fact, this was the method used in \cite{BvS}.)
In good cases one obtains a fourth order operator with MUM. In \cite{BK}
{\sc Batyrev} and {\sc Kreuzer} produced a list of promising candidate
Laurent polynomials. Some new Calabi--Yau operators were found in this
way.

The group of  {\sc Corti}, {\sc Coates}, and {\sc Kasprzyk} from
Imperial College in London has been persuing this approach on a larger scale,
systematically  using all reflexive  polytopes and Laurent polynomials with
special choice of the coefficients, \cite{CCGK}. From a preliminary run \cite{CCK}, $19$ new operators were found and one can reasonably expect many more to come from this approach.\\

{\em  Diagonals:} Not all Calabi--Yau operators arise from Laurent polynomials.
One obstruction comes from the fact that the numbers
\[a_n =[W^n]_{0}\]
coming from a Laurent polynomial satisfy {\em Dwork congruences}, \cite{SvS2} and \cite{MV}. The simplest of these imply that $a_n$ satisfy for each prime
number $p$ the congruence
\[ a_{n_0+n_1 p+\ldots n_k p^k}=a_{n_0}a_{n_1}\ldots a_{n_k} \mod p .\]

A more general concept is that of a {\em diagonal}.
If
\[ f=\sum a_{k_1 k_2 \ldots k_n}X_1^{k_1}X_2^{k_2}\ldots X_n^{k_n} \in \Q[[X_1,X_2,\ldots,X_n]]\]
is a power series in $n$ variables, then the {\em diagonal $\Delta_n(f)$ of $f$} is the power series in one variable obtained by only retaining the diagonal coefficients:
\[\Delta_n(f) := \sum_{k=0}^{\infty} a_{k k \ldots k} t^k \in \Q[[t]] .\]
As was shown by {\sc Christol} \cite{Chr}, the diagonals of rational functions
$P/Q, Q(0) \neq 0$ always satisfy a Fuchsian differential equation which is of geometric origin. More generally, by {\sc Lipschitz} \cite{Lip} the diagonal of any $D$-finite (i.e. holonomic) series is $D$-finite.

Note that if $W$ is a Laurent polynomial in $X_1,X_2,\ldots,X_n$, then its
constant term series is a diagonal of a rational function:

\[\sum_n [W^n]_0 t^n =\Delta_{n+1} \left(\frac{1}{1-X_0X_1\ldots X_n W(X_1,X_2,\ldots,X_n)}\right)\]

Any binomial sum can be converted into a representation as the diagonal of a rational function \cite{BLS}, hence the corresponding differential operator is always of geometrical origin in
the sense of \cite{And}.\\

In these constructions {\em integrality of the solution} is put in by construction, but the integrality of the mirror map and instanton numbers is for most operators {\em conjectural} and an experimental fact only.\\

{\bf \em Question 2:} How many Calabi--Yau operators do exist? Is their number  finite or infinite?\\

Of course, two operators that are related by a coordinate transformation or by multiplication with an algebraic function are to be considered as equivalent and we should count classes. Clearly, this is related to the question if there are finitely many or infinitely many distinct
topological types of Calabi--Yau threefolds, one of the big mysteries of the subject. It is not clear what to expect nor what to hope for.\\

A Calabi--Yau operator can be written in {\em $\Theta$-form} as
\[\mathcal{P}:= \Theta^4+tP_1(\Theta)+t^2 P_2(x)+\ldots+t^r P_r(\Theta),\]
where the $P_k$ are polynomials of degree four in $\Theta$ and we assume
$P_r \neq 0$. The number $r$ is then called the {\em degree} of the
Calabi--Yau operator $\mathcal{P}$. Over the last $12$ years {\sc Almkvist},
myself and others have been busy with collecting, simplifying and sorting
operators by degree, which is the simplest measure of complexity.
Our most recent list (August 2016) contains the the following operators:
\[
\begin{array}{|r|c|c|c|c|c|c|c|c|c|c|}
\hline
\textup{degree}         &1 &2 &3 &4 &5  & 6& 7& 8& 9&10\\
\hline
\textup{number of cases}&14&70&36&77&134&42&19&84&12&10\\
\hline
\end{array}
\]
\[
\begin{array}{|r|c|c|c|c|c|c|c|c|c|c|}
\hline
\textup{degree}         &11&12&13&14&15&16&17&18&19&20\\
\hline
\textup{number of cases}&20&17& 6& 9&0& 16& 0& 0& 0&1\\
\hline
\end{array}
\]
\[
\begin{array}{|r|c|c|c|c|c|c|c|c|}
\hline
\textup{degree}         &21&22&23&24&\ldots&32&\ldots&40\\
\hline
\textup{number of cases}&2&0&0&17&\ldots&1&\ldots&2\\
\hline
\end{array}
\]

Note that we are counting {\em operators together with the choice of a MUM-point}.
Some operators have more than one MUM-point, so these make, in transformed form,
multiple appearance on the list.  All listed cases are really different, as the
instanton numbers are different. However, it is conceivable that some of the operators 
of high degree are transformable to ones of lower degree, which would change the 
above table correspondingly. The operators are collected in a database that is accessible 
online, \cite{CYDB}.\\

{\bf \em Question 3:} For which operators do exist Calabi--Yau incarnations?\\

One might ask: does there exist an $A$-incarnation for a given operator
$\mP$? That is, does there exist a Calabi--Yau threefold $X$ with $h^{11}=1$ for which the instanton numbers are the instanton numbers of the operator?
\[ n_d(X)=n_d(\mathcal{P}) ?\]
For this we have only few examples and there are many operators which
the existence of such a manifold is in serious doubt.
For so-called {\em conifold operators} (see section 2.6) one can compute characteristic numbers like the Euler number from the monodromy of the operator, \cite{EvS}. However, there are a number of cases where this number turns out to be {\em positive}. See \cite{DeNu} for an example worked out in detail.

Does there exist a $B$-incarnation for $\mathcal{P}$? That is, does
there exist a  Calabi--Yau threefold $Y$ with $h^{12}=1$ for which
the Picard--Fuchs operator is $\mathcal{P}$? If this happens, we say that $\mathcal{P}$
has a strong $B$-incarnation (if $Y$ is even projective we call it a very strong
$B$-incarnation).  One could ask the operator to be a right factor of the Picard--Fuchs of a Calabi--Yau variety with $h^{12}>1$, which might be called a weak $B$-realisation. Differential operators having an integral solution are $G$-operators in the sense of \cite{And}; conjectures of {\sc Bombieri} and {\sc Dwork} say that such an operator is of geometric origin. For the cases where the coefficients $a_n$ have a representation as a binomial sum it is a theorem that they are of geometric origin.\\

Even if for some operators there do not exist strict $A$- or $B$-incarnations, it seems fruitful to consider each member of the list as describing something like a
{\em rank four Calabi--Yau motive over $\P^1$} and try to reconstruct as much as
possible of the geometry out of the differential operator alone.\\

\centerline{\bf \em Some recent examples}
\vskip 5pt
We report on two classes of examples of Calabi--Yau threefolds with
$h^{12}=1$ that are geometrically accessible and exhibit various
interesting phenomena. These examples will be discussed in two forthcoming
papers with {\sc Cynk} \cite{CvS3} and \cite{CvS4}. The computation of the
corresponding Picard--Fuchs operators became possible using the program of
{\sc Lairez}, \cite{Lai}. \\

{\bf Double octics.} A {\em double octic} is a threefold $Y$ that arises as
the double cover of $\P^3$ ramified over a surface $D \subset \P^3$ of degree 
eight.

If $D$ is smooth, $Y$ is a smooth Calabi--Yau threefold with
Hodge numbers $h^{11}=1$, $h^{12}=149$. If $D$ has singularities,
$Y$ is singular as well, but sometimes admits a crepant resolution.
Of particular interest is the case where $D$ is the union of eight planes:
as long as there are no fivefold points or fourfold lines in the 
configuration of planes, there exists a (projective) crepant resolution $\hat{Y}$ that can be obtained as covering of the blow-up of $\P^3$. By \cite{CvS0}, the infinitesimal
deformations of $\hat{Y}$ can be identified with the equisingular deformations of the
divisor $D$, which thus can be read of from the combinatorics of the intersection pattern
of the planes. 
In the thesis of {\sc  Meyer} \cite{Mey}, $63$ families of such double octics with $h^{12}(\hat{Y})=1$ were identified.
We determined for all these cases the corresponding Picard--Fuchs operator and some new operators were found this way. A particular beauty is the Calabi--Yau obtained from arrangement number $254$ of {\sc Meyer}. The octic $D$ is defined by the equation
\[xyzu(x+y+z+u)(u+y+tz)(zt+tu+x+y)(x+ty+zt) =0,\]
where $x,y,z,u$ are coordinates on $\P^3$ and $t$ is the parameter.
The Picard--Fuchs operator is a rather complicated operator of degree $12$ with Riemann symbol
\[\left\{
\begin{array}{ccccccccc}
0&1/2&1&\alpha&\beta&a&b&c&\infty\\
\hline
0&0&0&0&0&0&0&0&3/2\\
0&1&0&1&1&1&1&1&3/2\\
0&1&0&1&1&3&3&3&3/2\\
0&2&0&2&2&4&4&4&3/2\\
\end{array}
\right\} .
\]
The operator has {\em conifold points} at $1/2, \alpha, \beta$, where $\alpha,\beta$ are roots of $t^2-3t+1=0$, and
{\em apparent singularities} (so the local monodromy is trivial here) at $a,b,c$, roots of $8t^3-10 t^2+t-1=0$. At $0$, $1$ and at $\infty$ (after a quadratic
pullback) we have MUM points, with {\em three} different sets of instanton numbers.
\[
\begin{array}{|c|rrr|}
\hline
&0&1&\infty\\
\hline
n_1&288&128&4\\
\hline
n_2&59200&-4796&7/2\\
\hline
n_3&-8252768&341632&52\\
\hline
n_4&-1223488576&-31623118&500\\
\hline
n_5&585571467872&3395329408&2796\\
\hline
\end{array}
\]

On the A-side we expect {\em three} birationally distinct Calabi--Yau geometries with these (normalised) instanton numbers, which have equivalent derived categories!\\

In many other cases we obtain strong $B$-incarnations of operators that were known before. Also, there are many so-called {\em orphans} (see also section 2.5) and there are cases where the Picard--Fuchs operator is of order two. For details we refer to the forthcoming \cite{CvS3}.\\

It appears that there exist many operators in the list that have {\em two} points of maximal unipotent monodromy. In a recent series of papers \cite{HoTa1, HoTa2, HoTa3}, {\sc Hosono} and {\sc Takagi} described the beautiful geometry of the {\em Reye congruence Calabi--Yau threefold $X$}. The symmetric complete intersection of five divisors of degree $(1,1)$ in $\P^4 \times \P^4$ was considered in \cite{BvS}, but the Reye Calabi--Yau threefold $X$ arises from this complete intersection by dividing out the involution interchanging the $\P^4$-factors and has $h^{11}=1, h^{21}=26$. The corresponding Picard--Fuchs operator for the mirror family was described in \cite{BvS} and appeared as number $22$ in the AESZ-list \cite{AESZ}:
\[7^{2} \theta^4-\]
\[7 x\left(155\theta^4+286\theta^3+234\theta^2+91\theta+14\right)-\]
\[x^{2}\left(16105\theta^4+68044\theta^3+102261\theta^2+66094\theta+15736\right)+\]
\[2^{3} x^{3}\left(2625\theta^4+8589\theta^3+9071\theta^2+3759\theta+476\right)-\]
\[2^{4} x^{4}\left(465\theta^4+1266\theta^3+1439\theta^2+806\theta+184\right)+\]
\[2^{9} x^{5}\left((\theta+1)^4\right)\]
It has a second MUM-point at infinity, whose mirror appears to be the 
{\em double quintic symmetroid} $X'$, a double cover of the general linear 
symmetric determinant in $\P^4$.
So this is a second example quite similar to that of the Pfaffian and the Grassmannian Calabi--Yau and indeed in \cite{HoTa3} $X$ and $X'$ were shown to be derived equivalent.\\

{\bf Elliptic fibre products.} By blowing up the nine intersection points of a pencil of plane cubics we obtain a {\em rational elliptic surface} $\mathcal{E}$. By construction it admits a map $\pi: \mathcal{E} \lra \P^1$ and the fibres are identified with the cubics of the pencil; from the Euler number $\chi(\mathcal{E})=3+9$ we see that in general there will be $12$ nodal cubics in the family. As all cubics through eight of the base points also pass through the nineth, and
four of the eight points can be fixed, we see that the construction depends on 
eight parameters, and thus that there is {\em one condition} on the position of the $12$ singular fibres of a rational elliptic surface.
By specialisation of the construction, the singularities of these fibres may coalesce to form other Kodaira
types, but the sum of their Euler numbers will always add up to $12$. There are lists by {\sc Schmickler-Hirzebruch} \cite{Schm2} and by {\sc Herfurtner} \cite{Herf} that give all possible combinations of three and four Kodaira fibres; among them there are the six {\em Beauville surfaces} \cite{Beau1} with four fibres of type $I_n$.

In $1988$ {\sc Schoen} \cite{Scho} decribed a simple and very interesting class of Calabi--Yau threefolds by taking
the fibre product of two such rational elliptic surfaces $\mathcal{E}_i$, $i=1,2$:
\[ Y:=\mathcal{E}_1 \times_{\P^1} \mathcal{E}_2 \lra \P^1 .\]
If the sets of singular values $\Sigma_i \subset \P^1$ of $\mathcal{E}_i$ are disjoint, then $Y$ is a smooth
Calabi--Yau threefold that depends on $19=11+11-3$ parameters. The Euler number is equal to zero, as the fibres over a point of $\P^1$ all have Euler number zero, and indeed the Hodge numbers of $Y$ are $(h^{11},h^{12})=(19,19)$.
If, however, the fibrations $\mathcal{E}_i \lra \P^1$ have singular points in common, the threefold $Y$ aquires
singularities. For example, when an $I_n$ fibre meets an $I_m$ fibre, $Y$ aquires $n \cdot m$ singularities of type $A_1$. When we take a small resolution $\hat{Y}$ of these singularities, we obtain a smooth Calabi--Yau threefold whose
Hodge numbers can be determined easily from the singular fibres. In particular, there is a large number of cases where $h^{12}(\hat{Y})=1$, involving elliptic surfaces with up to six singular fibres.
In \cite{CvS1} we started exploring these examples, but it was only after {\sc Lairez}'s program \cite{Lai} became available that we were able to determine the most complicated of the corresponding Picard--Fuchs operators.\\

{\em Example:} We take for $\mathcal{E}_1$ the Beauville surface with fibres $I_6,I_3,I_2,I_1$ and as
$\mathcal{E}_2$ a surface with five singular fibres $I_8,I_1,I_1,I_1,I_1$ that depends on a single modulus-parameter $t$ (the Weierstrass equation for this surface is too complicated to write down here).
We identify the bases of the fibrations of $\mathcal{E}_1$ and $\mathcal{E}_2$
in such a way that three of the fibres of these
two families appear over the same point of $\P^1$ in the following way:
\[
\begin{array}{|c|cccccc|}
\hline
 &0&1&\infty&&&\\
\hline
\mathcal{E}_1&I_6&I_3&I_2&I_1&-&-\\
\hline
\mathcal{E}_2&I_8&I_1&I_1&-&I_1&I_1\\
\hline
\end{array}
\]
For generic choice of the modulus parameter $t$, the three ``free'' $I_1$ fibres will be disjoint and we obtain a Calabi--Yau threefold with $8\cdot 6+3\cdot 1+2\cdot 1=53$ $A_1$-singularities. A small resolution of these is a (non-projective) Calabi--Yau threefold $\hat{Y}$ with Euler characteristic $2\cdot 53=106$, so the Hodge numbers are $(h^{11}(\hat{Y}),h^{12}(\hat{Y}))=(52,1)$.
If the modulus-parameter $t$ is varied, the free $I_1$-fibres move, and for certain values further collisions of fibres do occur, leading to varieties with other types of singularities.
The Picard--Fuchs operator obtained is the most complex one encountered up to now. It made up number $500$ in the  list and was presented to {\sc Almkvist} on
occasion of his $3^4$ birthday, \cite{CvS2}. Currently, we are exploring the properties of the operators that can be obtained from this rich class of examples, \cite{CvS4}.

\section{\bf Some background}
In the second part of this paper we explain in some more detail certain of the concepts that were freely used in the first part.

\subsection{\em Differential equations of Fuchsian type}
We first briefly go over the basic properties of linear differential operators on the Riemann
sphere $\P^1=\C \cup \{\infty\}$ relevant for our discussion. We refer to {\sc Gray} \cite{Gra}
for an overview of the historical development of this very rich subject. The classical book
of {\sc Ince} \cite{Ince} contains a treasure of information and is still worth reading.

The set of singularities $\Sigma \subset \P^1$ of a differential operator
\[\mathcal{P}:=a_n(t) \frac{d^n}{dt^n}+a_{n-1}(t) \frac{d^{n-1}}{dt^{n-1}}\ldots+a_0(t) \in \C[t,\frac{d}{dt}]\]
are the zeros of $a_n(t)$ together, possibly, with the point $\infty$. In a neighbourhood of each
point $p \not\in \Sigma$ one find a basis of holomorphic solutions to the differential equation. A singular point $p \in \Sigma$ is said to be {\em regular singular} if all
solutions grow, at radial approach, at most as a power of the inverse distance to $p$.  {\sc Fuchs} \cite{Fuc} found a simple condition for this to happen: $0$ is a regular singular point of $\mathcal{P}$ if and only if $ord_0(a_i/a_n) \le n-i$.
 Differential equations with only regular singular points are called {\em regular singular} or {\em of Fuchsian type}. The solutions to such equations have an expansion of the form
\[ \sum a_{\alpha,k} t^\alpha\log^k(t)\]
that is convergent on a slit disc around each singular point.
The {\em Riemann symbol} of an operator summarises the information about
the local behaviour of the solutions near the singular points of a
differential operator. It consists of a table with columns indexed by the
singular points under which the corresponding  {\em exponents} are written. 
These exponents determine the local behaviour of the solutions at a singular
point and can be determined as follows.
If the operator $\mathcal{P}$ is written in $\Theta=t\frac{d}{dt}$-form
\[\mathcal{P}:=P_0(\Theta)+tP_1(\Theta)+t^2P_2(\Theta)+\ldots+t^d P_d(\Theta), \]
then the {exponents} of $\mathcal{P}$  at $0$ are just roots of the polynomial $P_0(\Theta)$.
The exponents of $\mathcal{P}$ at an arbitrary point $p$ are obtained by first
translating the point $p$ to the origin $0$ and writing the transformed operator
in $\Theta$-form and reading off the new $P_0$.
The exponents of $\mathcal{P}$ at $p=\infty$ are obtained using the reciprocal
transformation $t \mapsto 1/t$. When $\mathcal{P}$ is given
in $\Theta$-form, this is very easy operation, as it amounts to reversing the
sequence of polynomials $P_0,P_1,\ldots, P_d$ and replacing $\Theta$ by $-\Theta$, so that
the exponents of ${\mathcal P}$ at infinity are given by the negatives of the
roots of the polynomial $P_d$.\\

\subsection{\em Picard--Fuchs equations}
The differential equations that arise in algebraic geometry are usually those that are satisfied by integrals of rational or algebraic functions that depend on parameters and are called {\em Picard--Fuchs equations}. It was realised early on that such differential equations do have regular singularities. According to {\sc Houzel} \cite{Hou}, probably the first example was the equation found by {\sc Euler} \cite{Euler1}, that describes the
circumference of an ellipse with semi-axes of length $1$ and $\sqrt{1-t^2}$
as function of its eccentricity $t$:
\[ I(t):=4\int_0^1 \sqrt{\frac{1-t^2x^2}{1-x^2}}dx =2\pi \Big(1-\Big(\frac{1}{2} \Big)^2t^2-\Big(\frac{1\cdot 3}{2\cdot 4}\Big)^2\frac{t^4}{3}-\Big(\frac{1\cdot 3\cdot 5}{2\cdot 4\cdot 6}\Big)^2\frac{t^6}{5}-\ldots \Big) ,\]
which became later know as a {\em complete elliptic integral of the second kind}.
From the series expansion the differential equation satisfied by $I(t)$ is easily found to be
\[\mathcal{P}=\Theta^2-t^2 (\Theta-1)(\Theta+1)\]
which has
\[
\left\{
\begin{array}{cccc}
0 &1&-1&\infty\\
\hline
0&0&0&-1\\
0&1&1& 1
\end{array}
\right\}
\]
as Riemann symbol.
Another early example is the period of the mathematical pendulum ($L=g=1$)
with initial angle $\phi$ as function of $t:=\sin(\phi/2)$:
\[4 \int_0^1\frac{1}{\sqrt{(1-x^2)(1-t^2x^2)}}dx=
2\pi\Big(1+\Big(\frac{1}{2}\Big)^2t^2+\Big(\frac{1\cdot 3}{2\cdot 4}\Big)^2 t^4+\Big(\frac{1\cdot 3\cdot 5}{2\cdot 4\cdot 6}\Big)^2 t^6+\ldots \Big) ,\]
which is an elliptic integral of the first kind and satisfies the {\em Legendre
differential equation}
\[ \Theta^2-t^2(\Theta+1)^2\]
with Riemann symbol
\[
\left\{
\begin{array}{cccc}
0 &1&-1&\infty\\
\hline
0&0&0&1\\
0&0&0&1
\end{array}
\right\}   . \]

After initial work by {\sc Poincar\'e} and {\sc Picard} \cite{Pic}, it was with the work of {\sc Dwork} \cite{Dwo}
and {\sc Griffiths} \cite{Gri1}, \cite{Gri2} that integrals on higher dimensional manifolds were  studied systematically and methods to determine the Picard--Fuchs equation were developed.

The general geometric setting is most conveniently formulated as follows. One starts with a projective family
\[ f: \mathcal{Y} \lra \P^1 \]
and let $\Sigma \subset \P^1$ be the set of critical values of $f$, so that $f$
is a smooth map when restricted to $\P^1\setminus \Sigma$. For $t \in
\P^1 \setminus \Sigma$ the fibre $Y_t:=f^{-1}(t)$ is a smooth $d$-dimensional
variety; for $t \in \Sigma$ the fibre $Y_t$ will aquire a singularity.
Now choose for $p \in \P^1 \setminus \Sigma$ a $d$-cycle $\gamma_p \in H_d(Y_p)$.
Using the local topological triviality of $f$  over $\P^1 \setminus \Sigma$, we
can transport $\gamma_p$ to neighbouring fibres and obtain cycles $\gamma_t \in H_d(Y_t)$
for $t$ in a neighbourhood of $p$. If we choose a relative differential $d$-form
\[ \omega \in \Gamma(\mathcal{Y}, \Omega^d_{\mathcal{Y}/\mathbb{P}^1})\]
one can form the {\em period integral}
\[ \phi(t) = \int_{\gamma_t} \omega_{|Y_t} ,\]
which initially is defined in a neighbourhood of $p$, but which can, using the
local topological triviality of $f$, be extended along arbitrary paths in
$\P^1 \setminus \Sigma$. The finite dimensionality of the cohomology space $H^d_{dR}(Y_t)$
implies that $\phi$ satisfies a linear differential equation, called the {\em Picard--Fuchs equation}.
If the family $f:\mathcal{Y} \lra \P^1$ is given in terms of polynomial equations it is in principle
always possible to find this Picard--Fuchs equation, but it might not be simple to do so in practice.
For recent computer implementations see \cite{Mov} and \cite{Lai}.

\subsection{\em Local systems and monodromy}
\vskip 5pt
If $\mathcal{P}$ is a differential operator of order $n$ with set of singularities $\Sigma$, then the set of solutions to $\mathcal{P} \phi =0$ form, in the neighborhood of
each regular point $p$, a $\C$-vector space $\L_p$ of dimension equal to $n$: at $p$ a solution is uniquely determined  by the values of the first $n-1$ derivatives.
Hence, we obtain a {\em local system of solutions} $\L:=Sol(\mathcal{P})$ on $\P^1 \setminus \Sigma$.
If we choose a base point $p \in \P^1\setminus \Sigma$, the local system determines and is determined by a {\em representation of the fundamental group}
\[ T: \pi_1(\P^1\setminus \Sigma,p) \lra Aut(\L_p) \simeq Gl_n(\C),\;\gamma \mapsto T_{\gamma} \]
called the {\em monodromy representation} of $\mathcal{P}$,  which describes the behavour of
solutions of $\mathcal{P}$ under analytic continuation along closed paths. The image
of the fundamental group under $T$ is called the {\em monodromy group} of $\mathcal{P}$.
A choice of generators $\gamma_1,\gamma_2,\ldots, \gamma_r$ of $\pi_1(\P^1\setminus \Sigma,p)$ results in  a $r$-tuple of matrices $T_i:=T_{\gamma_i}$
\[ (T_1,T_2,T_3,\ldots,T_r) \in GL_n(\C)^r\]
that completely describes the local system. A change of base in $\L_p$ leads to
a simultaneous conjugation of all $T_i$.
These powerful ideas were introduced by {\sc Riemann} \cite{Rie} in his study
of the classical hypergeometric function $F(\alpha, \beta,\gamma; t)$. In fact
he determined the monodromy representation for the classical hypergeometric
differential operator and showed that it characterised the equation.
In his thesis, {\sc Levelt} \cite{Lev} found the monodromy representation of the higher hypergeometric functions. In fact, in this case the representation is uniquely determined by the Jordan type of the local monodromies around the singular points and thus represent the simplest examples of what is now called a {\em rigid local system}, \cite{NKatz}.\\

If the operator $\mathcal{P}$ is Fuchsian, the Zariski closure of the monodromy group is equal to the {\em differential Galois group} $Gal(\mathcal{P})$ that is introduced in the theory of {\sc Picard} and {\sc Vessiot} in analogy with the Galois group of an algebraic equation. We refer to the book
\cite{PutSin} for a detailed account.\\

In his famous 1900 ICM adress held in Paris, {\sc Hilbert} asked in his 21.~problem for the existence of a Fuchsian differential operator with prescribed monodromy representation. This is the so-called {\em Riemann--Hilbert problem} and is of central importance in contemporary mathematics; its solution is an interesting chapter in the history of mathematics, \cite{Ple1, Ple2, Beau2, AB, Del1}.\\

In the case of Picard--Fuchs operators, one starts with a projective family
\[ f: \mathcal{Y} \lra \P^1,\]
smooth over $\P^1 \setminus \Sigma$
and considers the direct image sheaf
\[ R^df_* \C_{\mathcal{Y}}\; .\]
It restricts to a local system $\H$ on $\P^1\setminus \Sigma$ and one has
\[ \H_t =H^d(Y_t,\C) .\]
In fact, the local system $\H$ has many further special properties. For instance,
we can consider the restriction $\H_{\Z}$ of $R^df_* \Z_{\mathcal{Y}}$ to $\P^1\setminus \Sigma$ and we have
\[ \H_{\Z} \otimes \C = \H\]
so $\H_{\Z}$ produces a lattice bundle inside $\H$ and so the monodromy
representation lands in $GL_n(\Z)$.
The behaviour of cycles under parallel transport was studied by {\sc Picard} \cite{Pic} and {\sc Lefschetz} \cite{Lef} and these works represent the first topological studies of higher dimensional algebraic varieties.
It led to the {\em Picard--Lefschetz} formula
\[ T_{\gamma} : H^d(Y_p) \lra H^d(Y_p), v \mapsto v \pm \langle v,\delta \rangle\delta\]
describing the cohomological monodromy that a cycle $v$ undergoes under
parallel transport along a path $\gamma$ that circumscribes (in the positive direction) a singular fibre that aquires a node ($A_1$-singularity) (see also \cite{Lam}, \cite{Looijenga}); $\delta$
is the {\em vanishing cycle}:  the class of a sphere that gets contracted when
passing to the singular fibre.
In general, all sorts of singularities might appear in the fibres. It is
a fundamental fact that the local monodromy transformation around any singular point $s$
\[ T_s : H^d(Y_t) \lra H^d(Y_t)\]
is always {\em quasi-unipotent}: there exist $m,k$ such that
\[ (T_s^m-I)^k=1 .\]
(In fact, one can take $k=d+1$.) This is called the {\em monodromy theorem}
and was first proven in \cite{Lan}. As a consequence, the exponents of
Picard--Fuchs operators are always rational.\\

\subsection{\em Self-duality}
The local systems $\H$ coming from geometry also have a build-in self-duality that  reflects {\em Poincar\'e duality} in the fibres: intersection of cycles in the fibres $Y_t$ defines a non-degenerate pairing
\[ \H_t \times \H_t \lra \C .\]
This leads to a self-duality property of the local system $\H$ and of the corresponding Picard--Fuchs operators $\mathcal{P}$.

Recall that the {\em adjoint} $\mathcal{P}^*$ of a differential operator is obtained by reading the operator backwards with alternating signs: if
\[\mathcal{P} =\sum_{i=0}^n a_{i}(t) \frac{d^i}{dt^i} \in \Q(t)\!\!\left[\frac{d}{dt}\right]\]
then
\[\mathcal{P}^*=\sum_{i=0}^n (-\frac{d}{ dt})^ia_{i}(t) \in \Q(t)\!\!\left[\frac{d}{dt}\right] .\]
This notion was introduced by {\sc Frobenius} in \cite{Fro}.

We say that an operator $\mathcal{P}$ is {\em essentially self-adjoint}
if there exists a function $\alpha \neq 0$ such that
\[ \mathcal{P} \alpha =\alpha \mathcal{P}^{*} .\]
If such $\alpha$ exists, it is easy to see that it has to satisfy the
differential equation
\[ \alpha' = -\frac{2}{n} a_1 \alpha .\]
So if the residues of $a_1$ are rational, $\alpha$  will be an algebraic function 
and the local system of solutions to an essentially self-adjoint operator
has a non-degenerate pairing with values in a rank one local system defined
by $\alpha$. It is symmetric if $n$ is odd and alternating if $n$ is even.

Only in the case that $\alpha$ is a rational function, we get an honest pairing and the differential Galois group $Gal(\mathcal{P})$ is contained in $Sp(n)$
($n$ even) or $SO(n)$ ($n$ odd).\\

For a fourth order differential operator
\[ \mathcal{P}:=\frac{d^4}{dt^4}+a_3(t) \frac{d^3}{dt^3}+a_2(t) \frac{d^2}{dt^2}+a_1(t) \frac{d}{dt} +a_0(t)  \in \Q(t)\left[\frac{d}{dt}\right] \]
the quantity
\[Q:=\frac{1}{2} a_2 a_3-a_1-\frac{1}{8}a_3^3+a_2'-\frac{3}{4}a_3 (a_3)'-\frac{1}{2} a_3''\]
was introduced in \cite{AZ} and taken as part of the definition of the
notion Calabi--Yau operator. It was shown in \cite{AZ} that the vanishing of
$Q$ is equivalent to the fact that the second order Wronskians of $\mathcal{P}$ satisfy an equation of order five rather than six and this is equivalent to $\mathcal{P}$ being essentially self-adjoint in the above sense.
As a consequence, the condition $Q=0$ does not always lead to operators with
$Gal(\mathcal{P}) \subset Sp(4)$, but only so after going to a cover defined
by the algebraic function $\alpha$  (which, in fact, is the unnormalised Yukawa coupling expressed in the original coordinate).
An example is the operator number $245$ from the AESZ-list \cite{AESZ}
\[ \mathcal{P}=
\Theta^4 -t (216\Theta^4+396\Theta^3+366\Theta^2+168\Theta+30)+ 36t^2 (3\Theta+2)^2(6\Theta+7)^2
\]
with instanton numbers
\[n_1=-6,\;\;\; n_2=-33,\;\;\;n_3=-170,\;\;\;n_4= -1029,\;\;\;n_5=-3246\]
for which $\alpha$ is
\[ \alpha(t)= \frac{1}{t^3(1-108t)^{11/6}} .\]
So only after going to a sixfold cover we do obtain an $Sp(4)$-operator.
We note that $\alpha$ is a rational function if and only if the exponents
at all singular points add up to an even integer, see e.g.\cite{Bog2}\\

\subsection{\em MUM and Hodge theory}
We say the operator $\mathcal{P}$ has a MUM point at $0$, if written in $\Theta$-form we have
\[\mathcal{P}=\Theta^4+tP_1(\Theta)+\ldots \]
In this case, the vector space $H_0$ of solutions on an (arbitrary small) slit
disc around the origin has a very special basis of solutions, called the
{\em Frobenius basis}
\[
\begin{array}{rcl}
y_0(t)&=&f_0(t)\\
y_1(t)&=&\log(t)y_0(t)+f_1(t)\\
y_2(t)&=&\frac{1}{2}\log(t)^2y_0(t)+\log(t)y_1(t)+f_2(t)\\
y_3(t)&=&\frac{1}{6}\log(t)^3y_0(t)+\frac{1}{2}\log(t)^2 y_1(t)+ \log(t) y_0(t)+ f_3(t), 
\end{array}
\]
where the $f_i$ are convergent power series with $f_0(0)=1$, $f_i(0)=0$, $i=1,2,3$.
We will also use the {\em scaled Frobenius basis}
\[ u_{3-k}:=y_k/(2\pi i)^k,\;\;\;k=0,1,2,3\]
The local monodromy around $0$ on the four dimensional vector space
\[H_0=\langle u_0,u_1,u_2,u_3\rangle\]
with respect to this scaled Frobenius basis is given by the matrix
\[ T_0:=\left(\begin{array}{cccc}
1&0&0&0\\
1&1&0&0\\
1/2&1&1&0\\
1/6&1/2&1&1
\end{array} \right) ,\]
which is unipotent with a Jordan block of maximal possible size,
hence {\em maximal unipotent monodromy}, that is for short, MUM.
%If $\mathcal{P}$ is self-dual, i.e. $Q=0$,
%one can show that the with respect to the $u$-basis the invariant symplectic fo%rm can be represented
%by the matrix
%\[ S:=\left(\begin{array}{cccc}
%0&0&0&1\\
%0&0&-1&0\\
%0&1&0&0\\
%-1&0&0&0
%\end{array} \right)\]

The local systems $\H$ on $\P^1 \setminus \Sigma$ that arise from algebraic
geometry have strong additional properties: the cohomology of the fibre $H^d(Y_t)$ carries a pure
Hodge structure of weight $d$: for each $t$ we have a Hodge decomposition
\[ H^d(Y_t)=\sum_{p+q=d}H^{p,q}_t,\;\; H^{p,q}_t=H^q(Y_t,\Omega^p)\]
that in fact depends nicely on $t$: the spaces of the {\em Hodge filtration}
\[ F_t^k =\sum_{p \ge k} H^{p,q}_t \]
form holomorphic vector bundles on $\P^1\setminus \Sigma$
\[ \mF^d \subset \mF^{d-1} \subset \ldots \subset \mF^0=\H \otimes \mO_{\P^1\setminus \Sigma} .
\]
One says that the local system $\H$ on $\P^1\setminus \Sigma$ underlies a
{\em variation of Hodge structures} (VHS). It is a fundamental fact
proven by {\sc Schmidt} \cite{Schm1} that one may extend this structure
defined on $\P^1 \setminus \Sigma$ over the punctures $s \in \Sigma$
to a  {\em mixed Hodge structure} (MHS):
for each $s \in \Sigma$ there is a $\Q$-vector space
$H_s=H_{\lim}^d(Y_s)$ of dimension equal to the rank of the local system $\H$. It can be defined as the sections of the $\Q$-local system over an arbitrary small
slit disc centered at $s$. Write the local monodromy as $T_s=U_sS_s$,
where $U_s$ is unipotent and $S_s$ is semi-simple, and define the
{\em monodromy logarithm} as
\[ N_s= -\log U_s=(1-U_s)+\frac{1}{2}(1-U_s)^2+\frac{1}{3}(1-U_s)^3+\ldots\]
The nilpotent endomorphism $N_s$ defines a {\em weight filtration}
\[ W_0 \subset W_1 \subset \ldots \subset W^{2d}=H_s, \]
which is characterised by the property that
\[ N^k_s: Gr^W_{d+k} \stackrel{\simeq}{\lra} Gr^W_{d-k} .\]
One can use the Hodge filtration $\mF^{\bullet}$ to define a limit Hodge
filtration $F^{\bullet}_s$ on $H_s$ and the fundamental theorem is that
for each $s\in \Sigma$ the triple $(H_s,W_{s,\bullet}, F^{\bullet}_s)$ is
a mixed Hodge structure: the filtration $F^{\bullet}_s$ defines a pure
Hodge structure of weight $k$ on the graded pieces $Gr^W_k H_s$.\\

In the geometrical case {\sc Steenbrink} \cite{Ste}
has constructed this mixed Hodge structure on $H_s$ using a semi-stable model
\[ \begin{array}{ccc}
D &\hookrightarrow & \mZ\\
\downarrow&&\downarrow\\
\{s\}& \hookrightarrow & \Delta\\
\end{array}
\]
over a disc $\Delta$. Here $D$ is a (reduced) normal crossing divisor
with components $D_i$. The complex of relative logarithmic differential forms
\[\Omega_{\mZ/\Delta}^{\bullet}(\log D)\]
can be used to describe the cohomology of the fibres and its extension to
$\Delta$. The complex comes with two
filtrations $F^{\bullet}$, $W_{\bullet}$, which induces filtrations on the
hypercohomology groups
\[ \H^d(\Omega_{\mZ/\Delta}^{\bullet}(\log D)\otimes \mO_D),\]
which then leads to the limiting mixed Hodge structure on $H_s$.
We refer to \cite{PetSte} for a detailed account.
Of particular relevance is the resulting {\em weight spectral sequence},
which expresses the graded pieces $Gr^W_k H_s$ in terms of the intersection
pattern of the exceptional divisors  and which degenerates at $E_2$.
The $E_1$-term is given by
\[
\begin{array}{ccc}
E_1^{p,q}&:=&\sum_k H^{q+2(p-k)}(D[2p-k])\\[3mm]
&=&H^{q+2p}(D[2p]) \oplus H^{q+2p-2}(D[2p-1]) \oplus  \ldots,
\end{array}
\]
where
\[D[k] := \coprod D_{i_0} \cap D_{i_1} \cap \ldots D_{i_k}\]
and where the sum runs over all indices $i_0 < i_1 <\ldots <i_k$.
In the diagram below, the stars indicate possible non-zero entries in the
$E_1$-page of the weight spectral sequence of a degeneration of a threefold.
\[\begin{array}{ccccccc}
*&*&*&*& & & \\
 &*&*&*& & & \\
 &*&*&*&*& & \\
 & &*&*&*& & \\
 & &*&*&*&*& \\
 & & &*&*&*& \\
 & & &*&*&*&*\\
\end{array}
\]
The differential runs horizontally, the operator $N$ acts on it and goes
two steps to the right and two steps down. There is a reflection
symmetry around the central point.\\

The bottom row $E^{p,0}_1$ can be identified with the complex
\[ 0 \lra H^0(D[0]) \lra H^0(D[1]) \lra H^0(D[2]) \lra \ldots \lra H^0(D[d]) \lra 0\]
where the differential is induced by the inclusion maps. The {\em dual intersection complex} $\Gamma$  has $0$-cells in bijection to the irreducible components of $D$, $1$-cells in bijection to the intersections of divisors, etc. So we see that the bottom row complex computes the cohomology of the dual intersection complex. Hence
\[ Gr^W_0 H^k_{\lim}(Y_0) = H^{k}(\Gamma) .\]

There are four possibilities for the mixed Hodge diamond of the
limiting mixed Hodge structures appearing for
variations of Hodge structures with $h^{3 0}=h^{2 1}=h^{1 2}=h^{0 3}=1$.
In the diagrams below the $k$-th row from the bottom gives the Hodge numbers of $Gr^W_k$; the operator $N$ acts in the vertical direction shifting downwards by two rows. The diagram is symmetric around the central vertical axis (by complex conjugation) and the central horizontal action (by symmetry of the  weight filtration). The numbers in each slope $=1$ (so SW-NE-direction) row of the diagram have to
add up to the corresponding Hodge number, so are all equal to $1$ in our case.
The cases that arise are:

{\em F-point}
\[
\begin{array}{ccccccc}
&&&0&&&\\
&&0&&0&&\\
&0&&0&&0&\\
1&&1&&1&&1\\
&0&&0&&0&\\
&&0&&0&&\\
&&&0&&&\\
\end{array}
\]
In this case $N=0$, so this happens if and only if the monodromy is of {\em finite order}. The
limiting mixed Hodge structure is in fact pure of weight three. This happens
in the mirror quintic at $\infty$, where the monodromy is of order five.\\

{\em C-point}
\[
\begin{array}{ccccccc}
&&&0&&&\\
&&0&&0&&\\
&0&&1&&0&\\
1&&0&&0&&1\\
&0&&1&&0&\\
&&0&&0&&\\
&&&0&&&\\
\end{array}
\]
In this case $N \neq 0$, $N^2=0$ and there is a single Jordan block.
The pure part $Gr^W_3$ is a rigid Hodge structure with Hodge numbers
$1,0,0,1$. Furthermore, $Gr^W_4$ and $Gr^W_2$ are one-dimensional and
are identified via $N$. This type appears when a Calabi--Yau threefold
aquires one or more ordinary double points, nowadays often called
{\em conifold points}, which explains our name $C$-type point for it.
In the mirror quintic this happens at $t=1/5^5$. But there are many
different kinds of singularties that lead to this mixed Hodge diamond.\\

{\em K-type point}

\[
\begin{array}{ccccccc}
&&&0&&&\\
&&0&&0&&\\
&1&&0&&1&\\
0&&0&&0&&0\\
&1&&0&&1&\\
&&0&&0&&\\
&&&0&&&\\
\end{array}
\]

In this case we also have $N \neq 0$, $N^2=0$ but there are two
Jordan blocks. In this case the pure part $Gr^W_3 =0$ and $Gr^W_4$,
$Gr^W_2$ are Hodge structures with Hodge numbers $1,0,1$, which
are identified via $N$. The Hodge structure looks like that of the
transcendental part of a K3-surface with maximal Picard number, which
explains our name $K$-point for it. This type of degeneration does
not appear in the family of the quintic mirror, but is common in other
examples. The holomorphic three-form  is destroyed by the singularties 
appearing in the fibre.\\

{\em MUM-point}
\[
\begin{array}{ccccccc}
&&&1&&&\\
&&0&&0&&\\
&0&&1&&0&\\
0&&0&&0&&0\\
&0&&1&&0&\\
&&0&&0&&\\
&&&1&&&\\
\end{array}
\]

Here $N^3 \neq 0$ and there is a single Jordan block of maximal size.
The Hodge structures $Gr^W_{2k}$ ($k=0,1,2,3$) are one-dimensional and
necessarily of Tate type. This happens for the quintic mirror at $t=0$
and is one of the main defining properties of Calabi--Yau operators.\\
So at a MUM-point, the resulting mixed Hodge structure is an iterated
extension of Tate--Hodge structures. {\sc Deligne} \cite{Del3} has shown
that the instanton numbers $n_1, n_2, n_3, \ldots$ can be seen to 
encode precisely certain {\em extension data} attached to the variation 
of Hodge structures near the MUM-point.\\

Now it is very well possible that in a family of Calabi--Yau threefolds no MUM-points appear. In \cite{Roh} first examples were given and in \cite{GG} a further
example was described, but the corresponding Picard--Fuchs equation was of second
order. In \cite{CvS1} an example with differential Galois group $Sp(4)$ was given
and, in fact, there are many more. {\sc Zudilin} \cite{Zud} proposed to call such operators {\em orphans}, as they do not have a MUM. For Calabi--Yau threefolds
appearing in one-parameter families without a MUM-point, it is not clear how to approach
the problem of constructing a mirror manifold, nor how to extract enumerative information
of the mirror manifold using the Picard--Fuchs equation. For this, one will need to understand the information hidden in the extension data near C-type and K-type points.\\

\subsection{\em Integrality properties}
So far most of the properties of the differential operator we discussed were purely
algebraic and rather easy to arrange for. For Calabi--Yau operators one
supplements these by further arithmetic  {\em integrality conditions}.
Initially in \cite{AZ} it was required that the operator has an integral
solution, but it is more natural to allow small denominators and ask for
$N$-integral solutions.

We will now explain in some detail the reasons for the integrality of the
normalised period near a MUM-point in the case it arises as Picard--Fuchs
operator of a family of Calabi--Yau varieties defined over $\Q$. It can be
seen as a generalisation of the celebrated {\em Theorem of Eisenstein}.
In the year $1852$ {\sc Eisenstein} \cite{Eis} reported at the meeting of the
{\em K\"oniglich Preu{\ss}ische Akademie der Wissenschaften zu Berlin}
on a curious general property of the power series development
of algebraic functions: if the power series
\[ \phi(t) =a_0+a_1 t + a_2 t^2+\ldots \in \Q[[t]] .\]
solves an equation $0=R(t,\phi(t))$ where $R \in \Z[x,y]$, then only
{\em finitely many primes appear in the denominators of the coefficients $a_i$}:
\[ \phi(t) \in \Z[\frac{1}{N}][[t]]\]
He mentions the example $\sqrt{1+t}$, where the replacement of $t$ by $4t$
turns the series into one with integral coefficients, a fact he considers
as well-known. On the other hand, in the series expansion of $\log(1+t)$
and $e^t$ any prime appears in a denominator of a coefficient, and the theorem of {\sc Eisenstein} implies the well-known fact that these series are not algebraic.\\

As usual, there is a prehistory that goes back to {\sc Euler.} In a letter  to {\sc Goldbach}, {\sc Euler} \cite{Euler2} reported on the counting of the number $A_n$ of ways to decompose an $n$-gon into triangles by drawing diagonals and found
\[
\begin{array}{|c|cccccccc|}
\hline
n&3&4&5&6&7&8&9&10\\
\hline
A_n&1&2&5&14&42&132&429&1430\\
\hline
\end{array}
\]
and gave
\[ \frac{1-2a -\sqrt{1-4a}}{2aa}\]
as the generating function for this sequence of numbers that nowadays are called
the {\em Catalan-numbers}.  In his reply, {\sc Goldbach} expressed his delight in the fact that this square root function apparently has  integral coefficients for its expansion. In the next letter, {\sc Euler} remarks that more generally the expansion of
\[ \sqrt[n]{1-n^2a}\]
in powers of $a$ has only integers as coefficients.\\

{\sc Eisenstein} did not write down a formal proof of his discovery,
but indicated that {\em once
the truth of the statement was recognised, it was easy to show its truth
by the method of undetermined coefficients}. What did he have in mind?\\

It was {\sc Heine} \cite{Hei} who gave a proof of a sharpened version of
{\sc Eisenstein}'s claim. {\sc Heine} also remarked that the series
\[1+ \frac{1}{3}t+\frac{1}{3^2}t^2+\frac{1}{3^9}t^3+\frac{1}{3^{16}}t^4+\ldots\]
was not excluded by the theorem of {\sc Eisenstein}, but nevertheless was not
algebraic.\\

{\bf Definition:} A series $\phi(t) \in \Q[[t]]$ is called {\em $N$-integral}
if there exist $c, N \in \N$ such that
\[ c\phi(Nt) \in \Z[[t]] .\]

{\sc Heine} proved the following statement:\\

{\bf Theorem of Eisenstein:} {\em Algebraic series are $N$-integral.}\\

Another proof was given by {\sc Hermite} \cite{Herm} and there is a very nice proof of the result using the theory of {\em diagonals} that we explain now.\\

Recall that the diagonal of
\[ f=\sum a_{k_1 k_2 \ldots k_n}x_1^{k_1}x_2^{k_2}\ldots x_n^{k_n} \in \Q[[x_1,x_2,\ldots,x_n]]\]
is the power series
\[\Delta_n(f) := \sum_{k=0}^{\infty} a_{k k \ldots k} t^k \in \Q[[t]] .\]
In this way we obtain a $\Q$-linear {\em diagonalisation map}
\[ \Delta_n: \Q[[x_1,x_2\ldots,x_n]] \lra \Q[[t]],\;\;f \mapsto \Delta_n(f)\]
We can consider the set of rational functions $R_n$
\[\frac{P(x_1,x_2,\ldots,x_n)}{Q(x_1,x_2,\ldots,x_n)},\;\;P,Q \in \Q [x_1,x_2,\ldots,x_n]\]
($\textup{with}\;Q(0) \neq 0$) that admit a power series expansion and
we say that a power series is an {\em $n$-diagonal}, if it is the diagonal
of such a rational function in $n$ variables, that is, if it belongs to
\[\Delta_n(R_n) \subset \Q[[t]]\;.\]

There is an obvious notion of $N$-integrality for series in
many variables. Rational functions (with rational coefficients) in
many variables are obviously $N$-integral: if we take $P$ and $Q$ with
integral coefficients, then we can take the denominator
of $P(0)/Q(0)$ as $N$. As diagonals of $N$-integral series
are clearly $N$-integral, we see that all $n$-diagonals are in fact
$N$-integral for some $N$.\\

{\bf Theorem:} (F\"urstenberg \cite{Fue})\footnote{The main interest of the paper \cite{Fue} lies, however, in the
statement that the situation is completely different over finite
fields: many more power series are algebraic, like
\[ \phi(t) =\sum_{n=0}^{\infty} t^{p^n},\]
which satisfies the equation
\[ \phi(t) = t+ \phi(t)^p , \]
hence is an algebraic series.
{\em If $K$ is a finite field, then all $n$-diagonals are algebraic.}}
{\em The 2-diagonals of rational functions are precisely the algebraic series.}

This theorem thus provides a natural proof of Eisensteins theorem.
Let us indicate the proof. If $F(x,y) \in R_2$, then one
can write
\[\phi(t):=\Delta_2(F)= \frac{1}{2 \pi i} \int_{\gamma} F(\zeta,\frac{t}{\zeta})
\frac{d \zeta}{ \zeta} .\]

The cycle $\gamma$ encloses some the poles of $F$ on the Riemann-surface
given by $xy=t$, so evaluating the integral by residues shows that $\phi(t)$ 
indeed is an algebraic function.

Conversely, if a series $\phi(t)$ solves $R(t,\phi(t))=0$,
where $R(x,y) \in \Z[x,y]$, $R(0,0)=0$, $\partial_y R(0,0) \neq 0$,
then it is a nice exercise to show that
\[\phi(t) = \Delta_2(F(x,y)) ,\]
where
\[F(x,y)=y^2 \frac{\partial_y R(xy,y)}{R(xy,y)} .\]

There is a generalisation of this result to more variables:\\

{\bf Theorem:} (Denef and Lipschitz  \cite{DL})
{\em The diagonal of algebraic power series in $n$ variables
is the diagonal of a rational function in $2n$ variables.}\\

From the proof of {\sc Eisenstein}'s theorem we see that the diagonalisation map has a natural
interpretation in terms of residues and integration over a vanishing cycle as was pointed out
in \cite{Del2}. Let us consider the following model situation:
$X:=\C^n$, $S=\C$ and the map
\[ p: X \lra S, (x_1,x_2,\ldots,x_n) \mapsto x_1x_2\ldots x_n=t\;.\]

The fibre $X_t:=p^{-1}(t)$, $t \neq 0$ is isomorphic to $(\C^*)^{n-1}$ and
contains the $(n-1)$-cycle $\Gamma_t$ defined by:
\[  |x_1|=t_1, |x_2|=t_2,\ldots,|x_n|=t_n  \subset X_t\]
where $t_1,t_2,\ldots,t_n$ are positive real numbers such that $t_1 t_2 \ldots t_n =|t|$.
There is a map
\[ \Omega_X^n \lra \omega_{X/S};\;\;\; \omega \mapsto Res {\Big(}\frac{\omega}{x_1x_2\ldots x_n - t}{\Big)} .\]
Now the statement is
\[\frac{1}{(2\pi i)^n}\int_{\Gamma_t} Res {\Big(} \frac{h dx_1 dx_2 \ldots dx_n}{x_1x_2\ldots x_n-t}{\Big)}=\Delta_n(h) .\]

We will now describe a general theorem that implies the $N$-integrality of the
invariant period for one-parameter families of Calabi--Yau manifolds near a MUM-point. The theorem has its roots in the work of {\sc Christol}, in
particular in the following example that can be found in \cite{Chr}.\\

The power series
\[ F(1/2,1/2,1;t)=1+\Big(\frac{1}{2}\Big)^2t+\Big(\frac{1\cdot 3}{2\cdot 4}\Big)^2t^2+\ldots\]
is the normalised period of the differential form
\[ \omega=Res\Big(\frac{dxdy}{f}\Big)=\frac{dx}{2y}\]
on the standard elliptic curve $E_t$ defined by
\[ f(x,y)=y^2-x(1-x)(x-t)=0\]
which for $t=0$ aquires a node.
The equation can be written in the form
\[t = x-\frac{y^2}{x(1-x)}=\frac{x^2(1-x)-y^2}{x(1-x)}=u \cdot v\]
where
\[u=\frac{x\sqrt{1-x}-y}{\sqrt{x(1-x)}},\;\;\;v=\frac{x\sqrt{1-x}+y}{\sqrt{x(1-x)}} .\]
Expressed in the coordinates $u$, $v$, the form  $\omega$ transforms to
\[ \sqrt{1-\Big(\frac{u+v}{2}\Big)^2} \frac{du}{u} \]
so we can represent the normalised period as the diagonal of an algebraic function
of two variables:
\[ \Delta_2 \left(\frac{1}{\sqrt{1-(\frac{u+v}{2})^2}}\right)=F(1/2,1/2,1;t)\]
and hence as a diagonal of a rational function of four variables.
In fact, in this case, it has even a representation as a diagonal of rational
function
\[\frac{4}{4-(x+y)(1+z)}\]
of three variables.\\

In \cite{And} (Theorem 2, p.185) this idea is generalised to higher dimensions.

Consider a projective family of $d$-dimensional varieties
\[ f: \mathcal{Y} \lra \P^1,\;\;Y_t=f^{-1}(t)\]
defined over $\Q$. We assume that $0$ is a MUM-point and that
$Gr^W_0 H_{lim}^d(Y_0)$ is one-dimensional. This implies that for $t \in \P^1\setminus \Sigma$ the Hodge space $H^{d,0}(Y_t)=H^0(\Omega_{Y_t}^d)$ is one dimensional. Pick a (rational) differential form
$\Omega \in H^0(\mY,\omega_{\mY}(*))$
and a cycle
\[ \gamma_t \in H_d(Y_t,\Z)\]
that is invariant under the monodromy.
We can form the {period function}
\[ \phi(t):= \int_{\gamma_t} \omega_t,\;\;\; \omega_t:=Res\Big(\frac{\Omega}{f-t}\Big)\]
that is defined in a sufficiently small disc around $0$. We can write
\[ \phi(t)=C y_0(t)\]
where the normalised period expands as
\[  y_0(t) :=1+a_1 t+a_2t^2+a_3t^3+\ldots \in  \Q[[t]] .\]

{\bf Theorem:} (Christol-Andr\'e)

{\em The series $y(t)$ is the diagonal of an algebraic function of $d+1$-variables.}\\

{\bf Proof:} We may take a semi-stable model
\[\mathcal{Z} \stackrel{\pi }{\lra} \mathcal{Y} \stackrel{f}{\lra} \Delta\]
of our family over a disc $\Delta$. We let $g:=f \circ \pi: \mZ \lra \Delta$ and
set $Z_t=g^{-1}(t)$. For $t \neq 0$ we have $Z_t \stackrel{\pi}{\approx} Y_t$.
The singular fibre $Y_0$ is replaced by a reduced normal crossing divisor
$D=\bigcup D_i \subset \mZ$.
We pull back $\omega$ to $\eta$ on $\mZ$ and can write the period function as
\[ \phi(t) =\int_{\delta_t} Res_{Z_t} \Big(\frac{\eta}{f-t}\Big) ,\]
where the cycle $\delta_t$ maps to $\gamma_t$ via $\pi$. From {\sc Steenbrink}s
construction of the limiting mixed Hodge structure on $H^d_{\lim}(Y_0)$ we obtain
from the weight spectral sequence a description $Gr^W_0 H^d_{\lim}(Y_0)$ in terms of the
intersections of the divisors $D_i$
\[ Gr^W_0 H^d_{\lim}(Y_0)=H^d(\Gamma) ,\]
where $\Gamma$ is the dual intersection complex of the divisors $D_i$.
By assumption, this space is one-dimensional so in particular, there must be points
$m \in \mathcal{Y}$ where $d+1$ divisors intersect.

For each such point, we may compare the behaviour of $f$ with the standard model.
Let $R=\mathcal{O}_{\mathcal{Y},m}$ the local ring of $\mathcal{Y}$ at one of these points $m$.
The equations defining the $d+1$ divisors meeting at $m$ may not belong to $R$, but in the Henselisation $\widehat{R}$ of $R$ we find elements $x_0,x_1,\ldots,x_{d}$ such that $D_i$ is  locally defined by $x_i=0$ and we can arrange that the map $f$ is given by
\[ (x_0,x_1,\ldots,x_n) \mapsto x_0x_1\ldots x_d=t\]
By considering the points with $|x_i|=t_i$ fixed, $t_0t_1\ldots t_d=|t|$, we obtain a real
$d$-dimensional torus $T_m(t) \subset \mathcal{Y}_t$ that vanishes at the point $m$ if $t \rightarrow 0$.
As the group $W_0H^d_{\lim}(Y_t)$ is supposed to be one-dimensional, all these tori $T_m(t)$
are homologous to a rational multiple of $\delta_t$.
Writing the form $\eta$ in terms of the coordinates $x_i$ we have
\[ \eta= h(x_0,\ldots,x_d) dx_0\wedge dx_1 \wedge \ldots \wedge dx_d\]
 with $h \in \widehat{R}$.
So
\[ \phi(t) =\int_{\delta_t} Res \Big( \frac{\eta}{g-t}\Big) =c\int_{T_m} Res \Big(\frac{h dx_0dx_1\ldots dx_d}{x_0x_1\ldots x_d-t } \Big)=c \Delta_{d+1} (h)\]

So the normalised period is indeed the diagonal of an algebraic function in $d+1$ variables. \hfill $\diamond$\\

Combining this with the theorem of {\sc Denef} and {\sc Lipschitz}  we get:\\

{\bf Corollary:} {\em The normalised monodromy invariant period near a MUM-point
of a Calabi--Yau $d$-fold is an  $2(d+1)$-diagonal, hence is $N$-integral.}\\

In particular, for Calabi--Yau threefolds, the period $y_0$ is an $8$-diagonal!\\

In the context of Calabi--Yau operators, one asks also for the integrality
of the mirror map $q(t)$. In some cases integrality of the mirror map has been shown by {\sc Lian} and {\sc Yau} \cite{LY}, {\sc Krattenthaler} \cite{Kra} and {\sc Delaygue} \cite{Dela} using purely number theoretic methods. But for the majority of cases, the integrality of the mirror map remains unproven.\\

The integrality of the $n_d$ is much deeper. In an A-incarnation, although supposed to count rational degree $d$ curves, $n_d$ is defined
in terms of Gromov--Witten invariants, so are a priori only in $\Q$. So the
integrality of the $n_d$, which was the biggest selling point of \cite{COGP},
is in the end the most mysterious aspect of the calculation. The conjectural
duality between Gromov--Witten theory and Donaldson--Thomas theory \cite{MNOP}
would  provide a natural explanation. The recent paper \cite{IP} provides a
proof of the integrality of the $n_d$ in case the operator has an
$A$-incarnation. It uses purely symplectic methods.\\

In the case of a B-incarnation one can use the link between the
action of Frobenius and the Yukawa coupling discovered  by {\sc Kontsevich, Schwarz} and {\sc Vologdsky}  \cite{KSV}, \cite{SchV}. A claim was made in \cite{Vol} that in the geometrical case the instanton numbers are $N$-integral, where the $N$ relates to primes of bad reduction in the semi-stable reduction.\\

In the thesis of {\sc Bogner} \cite{Bog} one finds the following
interesting operator:
\[
\mathcal{P}:=\Theta^4-8t(2\Theta+1)^2(5\Theta^2+5\Theta+2)+192
t^2(2\Theta+1)(3\Theta+2)(3\Theta+4)(2\Theta+3) .
\]
The operator has an integral solution
\[y_0(t)=1+16 t+576t^2+25600t^3+1220800t^4+\ldots,\]
integral mirror map
\[q=t+40t^2+1984t^3+106496t^4+\ldots,\]
integral Yukawa coupling
\[K(q)=1+8q-5632q^3-456064q^4-17708032q^5+\ldots\]
but the corresponding $n_d$'s are not integral:
$n_p$ has denominator $p^2$ for
\[p=3,5,7,11,13,17,19,\ldots \]
This is rather puzzling. The integrality of solution and mirror map clearly
indicate that we have a rank four Calabi--Yau motive and one would expect the 
general arguments for the integrality of \cite{Vol} to be applicable, but
apparently they are not. Maybe there is a different scaling of the coordinate
that repairs this defect, but up to now we have been unable to find it.\\

{\em \bf Questions:} \\

(i) Is there a proof of the integrality of the mirror map
along the same lines as the proof of integrality of the normalised period $y_0$? The higher dimensional strata in the divisor $D$ clearly will be relevant.\\

(ii) Constant terms series of Laurent polynomials are special cases of diagonals. The cases specially relevant to Calabi--Yau periods are the reflexive ones, and more generally those with a single interior point. Is there a similar theory of reflexive  diagonals?\\

(iii) For the constant term of the powers of a Laurent polynomial whose
Newton polyhedron contains a single interior point there are so-called
{\em Dwork congruences}, see \cite{SvS2} and \cite{MV}. Is there an analogue
for diagonals?\\

\subsection{\em Monodromy conjecture}

The monodromy group $\Gamma \subset Sp_4(\Z)$ appearing in one-parameter
families of Calabi--Yau threefolds is largely mysterious. The paper 
\cite{COGP} suggested that this group has infinite index in $Sp_4(\Z)$, but the arguments given were
insufficient. Only recently \cite{BT} this was proven to be the case for
seven of the $14$ hypergeometric cases. Somewhat surprisingly, in the
other seven cases the monodromy group turned out to have finite index,
\cite{SinV, Sin}. The situation for other one-parameter families
is under active study, \cite{HvS}. {\sc Bogner} and {\sc Reiter}
determined all symplectically rigid local systems of rank four \cite{BR};
strong B-realisations have been described recently in detail in \cite{DoMa}.\\

In the paper \cite{COGP} the explicit analytic continuation of solutions for the
operator were derived from the classical {\em Barnes integral representation}, which
is tied to the hypergeometric nature of the operator.
One of the striking features is the appearance of the number $\zeta(3)$ in combination
with the Euler number of the quintic. For this, the four loop correction to the sigma-model was
 made responsible. In that calculation, $\zeta(3)$ enters via the third derivative of the
$\Gamma$-function at $1$.\\
The paper \cite{HKTY} contained a general study of complete intersections in products of
projective spaces. The relevant holomorphic period function has nice expansions, whose coefficients
are expressible as a quotient of products of factorials, like the $a(n)=(5n)!/(n!)^5$ appearing
in the case of the mirror quintic. The solutions with log-terms can be obtained from the Frobenius
method, i.e. as the derivatives of $\sum a(n+\rho)x^{n+\rho}$ at $\rho=0$. If we replace all
factorials by $\Gamma$ using the relation $x!=\Gamma(1+x)$, we are naturally led to consider the
powerseries expansion
\[\frac{\Gamma(1+5x)}{\Gamma(1+x)^5} =1+\frac{5}{3}\pi^2x^2-40\zeta(3) x^3+\ldots\]
If we set
\[ h:= \frac{x}{2\pi i}\]
this can be written in the form
\[\frac{\Gamma(1+5x)}{\Gamma(1+x)^5} =\frac{1}{5}(5-\frac{50}{24}h^2-200\frac{\zeta(3)}{(2\pi i)^3} h^3+\ldots) .\]
Lo and behold, the coefficients of the expansion in $h$ contain the characteristic numbers of
the quintic $X$:
\[ \int_X H\cdot H \cdot H=5,\;\;\; \int_X c_1(X) \cdot H =0,\;\;\; \int_X c_2(X) \cdot H=50,\;\;\; \int_X c_3(X)=-200\;\;\; !\]
In the paper \cite{HKTY} these {\em remarkable identities} were observed to hold for all
complete intersections in products of projective spaces and a version of it was generalised to
the toric setting in \cite{HoLiYa}. Inspired by these facts and the formal similarity between
the Chern polynomial of the quintic
\[ \frac{(1+h)^5}{(1+5h)}\]
and the above series
\[ \frac{\Gamma(1+5x)}{\Gamma(1+x)^5}\]
{\sc Libgober} tried to find a general formulation of this relationship and introduced the Hirzebruch genus associated to the power series 
$\frac{1}{\Gamma(1+x)}$.\\

In \cite{EvS} we started computing monodromies numerically and discovered the systematic apearance of $\zeta(3)$ for general (fourth order)  Calabi--Yau operators. This led to the following general conjecture:\\.

\centerline{\bf Conjecture 1}
\vskip 5pt

{\em The monodromy matrices for a $Sp(4)$-Calabi--Yau operator with respect to
the scaled Frobenius basis $u_0,u_1,u_2,u_3$ have entries in
\[ \Q[\lambda] \]
where
\[\lambda:=\frac{\zeta(3)}{(2\pi i)^3} .\]}

Recall that the monodromy around $0$ on the solution space
\[H_0=\langle u_0,u_1,u_2,u_3\rangle\]
is always represented by the matrix
\[ T_0:=\left(\begin{array}{cccc}
1&0&0&0\\
1&1&0&0\\
1/2&1&1&0\\
1/6&1/2&1&1
\end{array} \right) .\]
Furthermore, one can show that the skew-symmetric form
\[ \langle -,-\rangle: H_0 \times H_0 \lra \C\]
determined by the conditions
\[ \langle u_0,u_3 \rangle=-\langle u_1,u_2\rangle=\langle u_2,u_1\rangle=-\langle u_3,u_0\rangle\]
is monodromy invariant.

A Calabi--Yau operator is called a {\em conifold operator} if the
exponents around the singular point $c$ nearest to $0$ are $0,1,1,2$ and
the monodromy around $c$ is a {\em symplectic reflection}.

So the monodromy $T_c$ around $c$ can be described in terms of a vector
 $S \in H_0$ via the formula
\[ T_c: H_0 \lra H_0;\;\; v \mapsto v-\langle v,S \rangle S .\]
We call $S$ the {\em reflection vector}.

For example, the operator of Candelas is a conifold operator with $c=1/5^5$
and the vector $S$ in the $u$-base is
\[  (5,0, \frac{25}{12}, -200 \lambda)^T\]
where we scale the skew-form by putting
\[\langle u_0,u_3 \rangle=\frac{1}{5} .\]
Now note that for the quintic Calabi--Yau $X \subset \P^4$ and $H \in H^2(X)$
the hyperplane class we have as mentioned before:
\[
\begin{array}{rcl}
5&=&\int_X H \cdot H \cdot H, \\
0&=& \int_X c_1(X) \cdot H\cdot H, \\
50&=&\int_X c_2(X) \cdot H, \\
-200&=&\int_X c_3(X) .
\end{array}
\]

One can verify for all cases where the operator arises from an $A$-incarnation
of a Calabi--Yau threefold $X$ with $h^{11}=1$, the operator indeed is a
conifold operator, and that the corresponding reflection vector is of the form
\[ (d,0, \frac{c}{24}, e \lambda)^T,\;\;\;\langle u_0,u_3 \rangle=\frac{1}{d}\]
and thus determines the characteristic numbers of $X$.
\[
\begin{array}{rcl}
d&=&\int_X H \cdot H \cdot H , \\
0&=&\int_X c_1(X)\cdot H \cdot H ,\\
c&=&\int_X c_2(X)\cdot H ,\\
e&=&\int_X c_3(X) .
\end{array}
\]

And even more, the reflection vectors of all further conifold operators
of the AESZ-list appear to be of this form, where $d,c,e$ are {\em integers}.
From this one may conjecture the existence of Calabi--Yau threefolds with
the given characteristic numbers. This was the main idea of the paper
\cite{EvS}. In the meantime, a few of these conjectured Calabi--Yau
varieties $X$ have indeed been found, but there is still a big gap.\\

There is a beautiful interpretation of the monodromy conjecture in terms of homological mirror symmetry supplemented by the $\Gamma$-class introduced by {\sc Libgober} \cite{Li} discussed above and developed further by {\sc Kontsevich, Katzarkov, Pantev} \cite{KKP} and {\sc Iritani} \cite{Iri}. We will sketch now that intriguing line of reasoning that, needless to say, is largely conjectural in general.

According to {\sc Kontsevich} \cite{Kon}, mirror symmetry should be understood as an equivalence of categories
\[ D(X) \stackrel{Mir}{\lra} F(Y)\]
and for the Calabi--Yau hypersurfaces $Y \subset \P^n$ a version of this has recently been proven by {\sc Sheridan} \cite{She}.

The category on the left is $D^b(Coh(X))$, the bounded derived category of coherent sheaves on $X$, on the right we have $D^{\pi}(Fuk(Y))$, the derived Fukaya category of $Y$. The objects in this category are represented by Lagrangian cycles (with local systems on them) in $Y$ and the category only depends on the
symplectic manifold underlying $Y_t$. The Hom-spaces in this
category are given by Floer homology $HF(L,L')$ groups; its Euler
characteristic is just the intersection product of the corresponding cycles:
\[\langle L,L' \rangle := [L]\cdot [L']\]
On the left hand side, the Euler characteristic of the Hom-spaces between
$\mE$ and $\mF$ in $D(X)$ is the Euler pairing
\[ \langle \mathcal{E},\mathcal{F} \rangle :=\sum_i (-1)^i \dim Ext^i(\mathcal{E},\mathcal{F}) .\]
Under mirror symmetry these should correspond:
\[ \langle \mE,\mF\rangle = \langle Mir(\mE),Mir(\mF) \rangle .\]

In the SYZ-picture of mirror symmetry, the spaces $X$ and $Y$ are related by $T$-duality: both $X$ and $Y$ are supposed to have the structure of dual three-torus fibrations over a common base.
In \cite{SYZ} it is argued that under $Mir$ the structure sheaf $\mO_p$ of a
point is  mapped to a lagrangian torus $\bf{T}$ (with local system on it) in $Y$, and the structure sheaf $\mO_X$ is mapped to a lagrangian sphere ${\bf S}$:
\[Mir(\mO_X)={\bf S},\;\;\;Mir(\mO_p)={\bf T}\]
and indeed
\[ \langle \mO_X,\mO_p \rangle = 1 =\langle \bf{S},\bf{T} \rangle .\]

But there is a certain asymmetry: if we map objects of $D(X)$
to $H^{ev}:=\bigoplus_k H^{2k}(X)$ via the Chern character, we can
express the Euler pairing as
\[ \int_X ch(\mathcal{E^*}) ch({\mathcal{F}}) Td(X) .\]
Objects of $F(Y)$ represented by lagrangian cycles map directly to $H^{odd}=H^3(Y)$ and the pairing $\langle L,L' \rangle $ is just given as an intersection number; no
tangential information like $Td(X)$ comes in.
To overcome this asymmetry, one modifies the Chern character by slipping in
a sort of square root of the Todd class.
Recall that the Todd class is the characteristic class coming from the power
series
\[\frac{x}{1-e^{-x}} .\]
The identiy $\Gamma(x)\Gamma(1-x)=\pi/\sin(\pi x)$ for the $\Gamma$-function
is equivalent to
\[ \Gamma(1+\frac{x}{2\pi i})\Gamma(1-\frac{x}{2\pi i})=e^{x/2} \frac{x}{1-e^{-x}} .\]
Now introduce the $\Gamma$-class as the characteristic class belonging to
power series expansion of
\[\Gamma(1+\frac{x}{2\pi i})=\exp(-\frac{\gamma}{2\pi i} x+\sum_{k=2}^{\infty} \frac{\zeta(k)}{(2\pi i)^k}\frac{x^k}{k})\]
So one puts:
\[\Gamma(T_X):=\prod_i \Gamma(1+\frac{\xi_i}{2\pi i})\]
where the $\xi_i$ are the chern roots of $T_X$.
Then one can write
\[\langle \mathcal{E},\mathcal{F} \rangle = ( \psi(\mathcal{E})^* \cdot  \psi(\mathcal{F})) ,  \]
where
\[ \psi(\mE):=\Gamma(T_X) \cup ch(\mF)\]
and the operation $*$ multiplies a component in $H^{2k}$ by $(-1)^k$.
So we are supposed to get a commutative diagram
\[
\begin{array}{ccc}
D(X)&\stackrel{Mir}{\lra}&F(Y)\\
\psi\downarrow&&\downarrow \phi\\
H^{ev}(X) &\stackrel{mir}{\lra}&H^{odd}(Y)
\end{array}
\]
where $mir$ is the cohomological mirror map. (The same diagram is discussed
in this context at various places in the literature, see e.g. \cite{Hos}).

The geometric monodromy of the family $Y_t$, $t \in \P^1\setminus \Sigma$,
can be realised symplectically and acts as auto-equivalences on $F(Y)$.
Via mirror symmetry there should be a corresponding action on $D(X)$.
And indeed, {\sc Kontevich} conjectured that the monodromy around the
MUM-point corresponds to the auto-equivalence
$\mO(H) \otimes: D(X) \lra D(X)$, whereas the {\em Seidel--Thomas twist} in the
spherical object $\mO_X$ corresponds to the symplectic Dehn twist along the
sphere $\bf{S}$.

But that implies that the reflection vector for the monodromy around the
conifold point should be equal to
\[S=mir(\psi(\mO_X))=\phi({ \bf S}) \in H^{odd}(Y) .\]
We will identify the space $H^{odd}$ with the solution space $H_0$,
spanned by the scaled Frobenius basis. Now we can work out everything!\\

In the basis $1,H,H^2,H^3$ for $H^{ev}(X)$ the operation of $\mO(H) \otimes$
is mapped via the Chern character to  multiplication with $e^h$, so is 
represented by exactly the same matrix as $T_0$ in the $u$-basis. 
Therefore, it is natural to put
\[ mir(h^k) = \alpha  u_k ,\]
where $\alpha$ is to be determined.
The element $\psi(\mO_X)=\Gamma(T_X)$ is computed to be
\[\psi(\mO_X)=1-\lambda_2 c_2(X)-\lambda_3 c_3(X) , \]
where
\[\lambda_k := \frac{\zeta(k)}{(2\pi i)^k} .\]
Writing this in terms of the basis $1,h,h^2,h^3$ and applying $mir$
we find for the reflection vector
\[S:=mir(\psi(\mO_X)) =\frac{\alpha}{d}(d u_0+\frac{c}{24} u_2-\lambda_3 e u_3) , \]
where
\[d=\int_X H \cdot H \cdot H,\;\;\;c=\int_X c_2(X) \cdot H,\;\;\;e=\int_X c_3(X) .\]
So we take $\alpha=d$.
Furthermore, $mir(\psi(\mO_p))$ represents the cohomology class of the
torus in $H^{odd}=H_0$. We find $\phi(\mO_p)=\frac{1}{d} H^3$, so that
\[T:=mir(\phi(\mO_p))=\frac{\alpha}{d}u_3=u_3 (=y_0) .\]
As we are supposed to have $\langle S,T \rangle=1$, we see that the right scaling
of the skew form indeed is obtained by putting
\[\langle u_0,u_3 \rangle=-\langle u_1,u_2 \rangle=\frac{1}{d}=\langle u_2,u_1\rangle=-\langle u_3,u_0 \rangle .\]

So, miraculously, everything fits \footnote{This I realised after a talk by Kontsevich (Vienna, 2008), where he explained the $\Gamma$-class.}
and completely explains the  structure of the conifold reflection
vector in the Frobenius basis. Usually a Calabi--Yau operator has more
conifold points and the corresponding reflection vectors should arise 
from other spherical objects in $D(X)$. As before, we obtain vectors of 
the same shape
\[ \frac{1}{d}(d,a,\frac{c}{24},e \lambda_3)^T\]
but now usually $a \neq 0$. For more complicated monodromy transformations
we do not have a real argument, but there is little reason to doubt the
general principle.\\

If the monodromy is not in $Sp(4)$, there usually appear algebraic
numbers and the monodromy matrices appear to be contained in
\[\overline{\Q}[\lambda] .\]
As an example, take the operator number 245 from the AESZ-list \cite{AESZ}
mentioned earlier:
\[\Theta^{4}-t(216\theta^{4}+396\Theta^{3}+366\Theta^{2}+168\Theta+30)+
36 \left( 3\,\Theta+2 \right) ^{2} \left( 6\,\Theta+7 \right) ^{2}\]
with Riemann symbol
\[\left\{
\begin{array}{ccc}
0&1/108&\infty\\
\hline
0&0&2/3\\
0&1/6&2/3\\
0&1&7/6\\
0&7/6&7/6\\
\end{array}
\right\} .
\]

For the matrix of the monodromy around $1/108$ in the Frobenius basis
we find
\[
\frac{\sqrt{3}+i}{2}
\left( \begin {array}{cccc} 1/2\,\sqrt {3}-72\,\sqrt {3}\lambda &-\sqrt {3}&0
&18\,\sqrt {3}\\\noalign{\medskip}-1/6\,\sqrt {3}&1/2\,\sqrt {3}&-3\,
\sqrt {3}&0\\\noalign{\medskip}4\,\sqrt {3} \lambda &1/12\,\sqrt {3}&1/2\,
\sqrt {3}&-\sqrt {3}\\\noalign{\medskip}-{\frac {1}{72}}\, \left( 1+
20736\,{\lambda}^{2} \right) \sqrt {3}&-4\,\sqrt {3}\lambda&-1/6\,\sqrt {3}&1/2\,
\sqrt {3}+72\,\sqrt {3}\lambda\end {array} \right) , \]
which indeed has order six.

One can generalise this to higher order operators:\\

\centerline{\bf Conjecture 2}
{\em
The monodromy matrices for a Calabi--Yau operator of order $n+1$
with respect to the scaled Frobenius basis have entries in
\[ \overline{\Q}[\lambda_2,\lambda_3,\ldots,\lambda_n] ,\]
where
\[\lambda_k:=\frac{\zeta(k)}{(2\pi i)^k} .\]
}
This conjecture has been verified numerically for almost all
Calabi--Yau operators of the list. For those cases that can be related to
specialisations of hypergeometric functions of more variables one can
in principle {\em prove} these numerical results. For details we refer to
the thesis of {\sc Hofmann}, \cite{Hof}.\\

It turns out that it is natural to formulate conjectures about the $\Gamma$-class in the context of Fano manifolds, as was done in the beautiful paper \cite{GGI}. The so-called {\em $\Gamma$-conjectures} formulated there can be motivated from mirror symmetry \cite{GI} and they imply the above monodromy conjecture in those cases where the Calabi--Yau manifold has an $A$-incarnation as a complete intersection in a Fano-manifold for which the $\Gamma$-conjectures are proven. In the paper \cite{GZ} the $\Gamma$-conjectures were verified for the case of Fano threefolds of Picard rank equal to one.\\

In this overview paper we have touched upon various aspects of
Calabi--Yau operators. We had to leave out several important topics, most notably be $p$-adic story, which involves the Dwork congruences \cite{SvS2, MV}, the computations of the local L-factors \cite{SvS1, Sam} and the  $p$-adic analogue of the $\Gamma$-conjectures that lead to the appearance of the $p$-adic analogue of $\zeta(3)$. Also higher genus instanton numbers for Calabi--Yau operators were completely left out of this account.\\

{\bf Acknowledgement:} This paper is an extended write-up of a talk given
at the conference {\em Uniformization, Riemann--Hilbert Correspondence, Calabi--Yau Manifolds, and Picard--Fuchs Equations} held at Mittag-Leffler july 6-july 10, 2015. I thank the organisers for inviting me to this nice conference and giving me an opportunity to speak. I also thank the other participants at the conference for their presentations. Furthermore, thanks to my collegues {\sc Gert Almkvist}, {\sc Slawek Cynk}, {\sc Vasily Golyshev} and {\sc Wadim Zudilin} for constant inspiration and my former students {\sc Michael Bogner}, {\sc J\"org Hofmann}, {\sc Christian Meyer} and {\sc Kira Samol} for their courage to embark on uncertain projects related to Calabi--Yau operators and bringing them to a good end. 
I also thank {\sc Pierre Lairez} bringing \cite{Hou} and \cite{Euler1} to my 
attention and for writing a very useful computer program. Many thanks also to 
the unknown referee for pointing out certain essential references that I 
missed in an earlier version of this text.

\end{document}